\documentstyle{article}
\input psfig.sty  
\newtheorem{definition}{Definition}[section]

\newtheorem{theorem}[definition]{Theorem}
\newtheorem{proposition}[definition]{Proposition}

\newtheorem{remark}[definition]{Remark}
\newtheorem{example}[definition]{Example}

\newcommand{\sect}[1]{\section{#1}\setcounter{equation}{0}}

\def \trait (#1) (#2) (#3){\vrule width #1pt height #2pt depth #3pt}
\def \qed{\hfill
        \trait (0.1) (6) (0)
        \trait (6) (0.1) (0)
        \kern-6pt
        \trait (6) (6) (-5.9)
        \trait (0.1) (6) (0)
\medskip}

\newcommand{\dx}{\,dx}
\newcommand{\dt}{\,dt}
\newcommand{\dy}{\,dy}

\newcommand{\ie}{; {\it i.e., }}
\newcommand{\e}{\varepsilon}

\newcommand{\lp}{{\rm L}^p}

\newcommand{\wup}{{\rm W}^{1,p}}

\newcommand{\wuu}{{\rm W}^{1,1}}

\newcommand{\orm}{(\Omega;\R^m)}
\newcommand{\loc}{{\rm loc}}

\newcommand{\Hom}{{\rm hom}}     %%%% = vecchio \hom

\newcommand{\xe}{{x\over \e}}

\newcommand{\ave}{-\hskip -.38cm\int}

\newcommand{\Om}{\Omega}

\font\tenmsb=msbm10
\font\sevenmsb=msbm7
\font\fivemsb=msbm5
\newfam\msbfam
\textfont\msbfam=\tenmsb
\scriptfont\msbfam=\sevenmsb
\scriptscriptfont\msbfam=\fivemsb
\def\Bbb#1{{\fam\msbfam\relax#1}}

\def\R{\Bbb R}
\def\M{\Bbb M}
\newcommand{\rr}{\Bbb R}

\newcommand{\Rm}{\R^m}

\newcommand{\NN}{{\Bbb N}}

\newcommand{\ZZ}{{\Bbb Z}}

\newcommand{\MM}{{\Bbb M}^{m\times n}}

\newcommand{\Aa}{{\cal A}}           %%%% = vecchio \AA
\newcommand{\Ao}{\Aa(\Omega)}

\def\D{D_{\alpha}}
\def\x{x_{\alpha}}

\def\F{\overline F}
\def\to{\rightarrow}

\title{Homogenization of oscillating boundaries and applications to thin
films}

\author{
{\sc Nadia Ansini} and {\sc Andrea Braides}
\\ {\normalsize SISSA, via Beirut 4, 34013 Trieste, Italy}
\\ {\normalsize Email: ansini@sissa.it, braides@sissa.it}
}

\date{}

\begin{document}

\maketitle

\sect{Introduction}
The study carried on in this paper draws its motivation from the problem of
the asymptotic
description of nonlinearly elastic thin films with a fast-oscillating profile.
The behaviour of such films is governed by an elastic
energy, where two parameters intervene: a first parameter $\e$ represents
the thickness of the thin film and a second one $\delta$ the scale of
the oscillations.
The analytic description of the elastic energy is given by a functional of
the form
\begin{equation}\label{1.01}
E_{\e,\delta}(u)=\int_{\Omega(\e,\delta)} W(Du)\dx,
\end{equation}
where the set $\Omega(\e,\delta)$ is of the form
\begin{equation}\label{1.02}
\Omega(\e,\delta)=
\Bigl\{
x\in\R^3:\ |x_3|< \e\,f\Bigl({x_1\over\delta},{x_2\over\delta}\Bigr),
(x_1,x_2)\in\omega
\Bigr\},
\end{equation}
with $f$ is a bounded $1$-periodic function which parameterizes the boundary
of the thin film, which then has periodicity $\delta$.
It is convenient to scale these energies by a change of variables
and consider the functionals
\begin{equation}\label{1.1}
E_{\e}^\delta(u)=\int_{\Omega(\delta)}
W\Bigl(D_1u,D_2u,{1\over\e}D_3u\Bigr)\dx,
\end{equation}
where now
\begin{equation}\label{1.2}
\Omega(\delta)=
\Bigl\{
x\in\R^3:\ |x_3|< f\Bigl({x_1\over\delta},{x_2\over\delta}\Bigr),
(x_1,x_2)\in\omega
\Bigr\}.
\end{equation}
In this way we separate the effects of the two parameters $\e$ and $\delta$.

In a recent paper by Braides, Fonseca and Francfort \cite{BFF}
a general compactness result for functional of thin-film type has been
proven which comprises energies of the form (\ref{1.1}), showing that, with
fixed
$\delta=\delta(\e)$, upon possibly extracting a subsequence, the family
$E_\e^{\delta(\e)}$ converges in the sense of De Giorgi's $\Gamma$-convergence
as $\e\to 0$ to a $2$-dimensional energy, which, if $\delta(\e)\to 0$ as
$\e\to 0$,
can be identified with a $2d$-functional of the form
\begin{equation}\label{1.3}
E(u)=\int_{\omega} \widetilde W(D_1u,D_2u)\dx.
\end{equation}
In many cases it is possible to describe $\widetilde W$ explicitly in terms
of $W$ and $f$, and as a consequence to prove that no passage to a
subsequence is
necessary.
When $f=C$ is constant ({\it i.e.}, the profile of the thin film is flat,
and hence there is no real dependence on $\delta$) the description
of the energy density $\widetilde W$ has been given by Le Dret and Raoult
\cite{LDR}
who proved that $\widetilde W= 2C\,Q_2\overline W$; here $Q_2$ denotes
the operation of
$2d$-quasiconvexification, and $\overline W$ is obtained from $W$ by
minimizing in the
third component. An equivalent formula, of `homogenization type', is given
in \cite{BFF}
(see also \cite{BD}). If $\overline W \neq Q_2 \overline W$
({\it i.e.}, $\overline W$ is not quasiconvex)
then both formulas underline the formation of microstructures generated
by the passage to the limit.
When $f$ is not constant, then the function
$\widetilde W$ depends
on the behaviour of $\delta$ with respect to $\e$. The case when
$\delta=\e$ (or more
in general when
$\delta/\e$ converges to a constant) has been treated in \cite{BFF},
where it is shown that a homogenization type formula for $\widetilde W$
can be given. The same
method can be used
when $\delta>>\e$; in this case the recipe to obtain $\widetilde W$ is the
following:
first, keep $\delta$ fixed and apply the Le Dret and Raoult procedure,
considering
the thickness of the thin film as a parameter. The output of this
procedure is a
$2$-dimensional
energy of the form
\begin{equation}\label{1.4}
E^\delta(u)=\int_{\omega} 2f\Bigl({x_1\over\delta},{x_2\over\delta}\Bigr)
Q_2 \overline W(D_1u,D_2u)\dx.
\end{equation}
We can then let $\delta$ tend to $0$, and apply well-known homogenization
procedures
(see \cite{BD}) obtaining a limit functional, which turns out to be the
desired one.
In the case $\delta<<\e$ it is possible to make an {\it ansatz} in the same
spirit, arguing
that the limit $E$ can be obtained in the following two steps:

(1) (Homogenization of sets with oscillating boundaries)
First consider $\e$ as fixed, and let $\delta\to 0$, to obtain
a limit functional of the form
$$
E_\e(u)=
\int_{\omega\times (-1,1)}W_\Hom\Bigl(x_3,D_1u,D_2u,{1\over\e}D_3u\Bigr)\dx
$$
(we consider the normalized case $\sup f=1$).

Note that in this case an additional dependence on $x_3$
is introduced, which may underline a loss of coerciveness of the function
$W_\Hom$
for certain values of $x_3$. The form of $W_\Hom$ will depend on
$W$ and on the sublevel sets of $f$;

(2) (Thin film limit) Let $\e\to 0$ and generalize the method of \cite{BFF}
to non-coercive functionals. In this way we obtain a limit energy density
\begin{eqnarray*}
\overline W_{\Hom}(\overline F)&=&\inf_{k\in\NN}\inf
\Bigl\{{1\over k^2}\int_{(0,k)^2\times(0,1)}
W_{\Hom}(x_3, Du+(\overline F,0))\dx:
\\ \nonumber&&\qquad\qquad
u\in\wup_\loc((0,1)^3;\R^3),\
u\ k\hbox{-periodic\, in}\, (x_1,x_2)\Bigr\}\,.
\end{eqnarray*}
Note that the dependence on $x_3$ implies that
the simpler method of \cite{LDR} cannot be applied to this situation.

A partial result in this case has been obtained by Kohn and Vogelius
\cite{KV} who dealt with linear operators.

\bigskip
The purpose of this work is twofold. First, we give a general theory for the
homogenization of non-convex energies defined on sets with oscillating boundaries
by generalizing the application of the direct methods of $\Gamma$-convergence
to homogenization as described in \cite{BD}. We clarify and prove statement (1) 
above, by showing that the functionals $E_\e$ are defined on a `degenerate 
Sobolev Space' that can be described by proving an auxiliary
convex-homogenization result. The formula for $W_\Hom$ can be obtained by
solving a possibly degenerate localized $3d$-homogenization problem. 
In the case of convex $W$ the determination of $W_{\Hom}(t,\overline 
F)$ for fixed $t\in(-1,1)$  essentially amount 
to solving a $2d$-homogenization problem with an 
energy which is coercive only on the set $E_{t}=\{(x_1,x_{2})\in \rr^{2}:
f(x_1,x_{2})>|t|\}$, while in the general non-convex case the problem defining
$W_{\Hom}(t,\overline F)$ is genuinely three dimensional.
We state and prove these results in a general $n$-dimensional setting
(for some related problems in the convex setting see 
{\it e.g.} \cite{BC}). 

The second goal of the paper is to prove that by following steps (1) and (2)
above we indeed obtain the description of $\widetilde W$. Even though 
this is an intrinsically vectorial problem, and hence the `natural' 
structural condition on $W$ is quasiconvexity, we have been able to prove this 
result only with the additional hypothesis that $W$ is convex. The 
technical point where this assumption is needed is the separation of 
scales argument, which assures that, essentially, homogenization comes 
first, followed by the thin film $3d$--$2d$ limit. In general
problems where only  quasiconvexity is assumed
 this point is usually proved by a compactness argument
which uses some
equi-integrability properties of gradients of optimal sequences for the 
homogenization derived from the growth conditions on the energy 
density (see {\it e.g.} Fonseca M\"uller Pedregal \cite{FMP}; 
for the use of this argument in the 
framework of iterated homogenization see \cite{BD} Chapter 22;
for an application to heterogeneous thin films with flat profile see 
Shu \cite{Shu}).
In the case of thin films with fast-oscillating profiles, this technique 
cannot be used since we have a control on the gradients of optimal sequences 
only on varying wildly oscillating domains. In the convex case though, optimal 
sequences for the homogenization can be obtained simply by scaling one single 
periodic function, and hence their gradients automatically enjoy 
equi-integrability properties. Note that this difficulty is similar 
to those encountered when dealing with higher-order theories of thin 
films. In that case the necessary compactness properties can be obtained by 
adding a small perturbation with higher-order derivatives
(as in the paper by Bhattacharya and James \cite{BJ}).
We do not follow this type of argument since even a singular 
perturbation by higher-order gradients might interact with the 
homogenization process, as shown by Francfort and M\"uller \cite{FrM}.
%%% in the case when $W$ is {\it convex}. Under this additional 
%%% hypothesis it is well known that the homogenization formula can be stated
%%% by using a single auxiliary cell problems; this fact allows us to avoid
%%% compactness arguments, which seem to be a difficult 
%%%issue for families of functions
%%% defined on sets with highly oscillating profile. 
More applications of $\Gamma$-convergence arguments to thin films theory can be 
found in \cite{BB, BF}.

\sect{Notation and Preliminaries}
\setcounter{equation}{0}

In the sequel, $n,m\in \NN$ with $n\ge 2$, $m\ge 1$.
If $x \in \R^n$ then $\x = (x_1,\ldots , x_{n-1}) \in \R^{n-1}$
is the vector of the first $n-1$ components of $x$, and $\D =
\left(\frac{\partial}{\partial x_1},\ldots ,
\frac{\partial}{\partial x_{n-1}}\right)$.
If $\Omega$ is a open subset of $\R^n$
we denote by $\Ao$ the family of all open subsets of $\Omega$.

The notation $\MM$ stands for the space of $m\times n$ matrices.
Given a matrix $F \in \MM$, and following the notation introduced in 
\cite{LDR}, 
we write $F = (\overline {F} | F_n)$, where $F_i$ denotes the $i$-th column
of $F$, $1 \leq i \leq n$, and $\overline {F} = (F_1,\ldots , F_{n-1}) 
\in \M^{m\times n-1}$ is the matrix of the first $n-1$ columns of $F$.
$\overline F$ denotes also $(\overline F,0)$ when no confusion arises.

The {\it characteristic function} of a set $E\subset\R^n$ is denoted
by $\chi_E$, and the $N$-dimensional Lebesgue measure in $\R^N$ is
designated as ${\cal L}^N$. We use standard notation for Lebesgue and 
Sobolev spaces. The letter $c$ will stand for an arbitrary
fixed strictly-positive constant.

\smallskip
We recall the definition of De Giorgi's {\it $\Gamma$-convergence}
in $\lp$ spaces, $1  \leq p < +\infty$.  Given a family of
functionals
$J_j:\lp(\Omega;\Bbb R^m) \to [0,+\infty)$,
$j\in\NN$,
 for $u\in \lp(\Omega;\Bbb R^m)$ we define
$$
\Gamma\hbox{-}\liminf_{j\to+\infty}
J_j(u)=\inf\Bigl\{\liminf_{j\to+\infty} J_j(u_j)\ :\ u_j\to u
\hbox{ in }\lp\orm\Bigr\},
$$ and
$$
\Gamma\hbox{-}\limsup_{j\to+\infty}
J_j(u)=\inf\Bigl\{\limsup_{j\to+\infty} J_j(u_j)\ :u_j\to u
\hbox{ in }\lp\orm\Bigr\}.
$$ 
If these two quantities coincide then their common value is
called the
$\Gamma$-{\it limit} of the sequence $(J_j)$ at $u$, and is denoted
by
$\Gamma\hbox{-}\lim_{j\to+\infty} J_j(u)$. It is easy to check
that $\l =\Gamma\hbox{-}\lim_{j\to+\infty} J_j(u)$ if and only
if

(a) for every sequence $(u_j)$ converging to $u$ in
$\lp(\Omega;\R^m)$ we have
$$
\l \le\liminf_{j\to+\infty} J_j(u_j);
$$

(b) there exists a sequence $(u_j)$ converging to $u$ in
$\lp(\Omega;\R^m)$ such that
$$ \l
\ge\limsup_{j\to+\infty} J_j(u_j).
$$ We say that $(J_\e)$ {\it $\Gamma$-converges to $\l$  at $u$ as
$\e\to 0$} if for every sequence of positive numbers $(\e_j)$
converging to $0$ there exists a subsequence $(\e_{j_k})$ for which
$$ 
l=\Gamma\hbox{-}\lim_{k\to+\infty} J_{\e_{j_k}}(u).
$$ 
We recall that the $\Gamma$-upper and lower limits defined above are
$\lp$-lower semicontinuous functions.

For a comprehensive study of
$\Gamma$-convergence we refer to the book of Dal Maso \cite{DM}
(for a simplified introduction see \cite{B}), 
while a detailed analysis of some of its applications to
homogenization theory can be found in \cite{BD}.

\sect{The direct method of $\Gamma$-convergence}
In the sequel we will repeatedly apply some variants of the so-called
direct method of
$\Gamma$-convergence to homogenization problems, which consists in combining
localization and integral representation procedures to obtain compactness
theorem for
classes of integral functional. This method in the version which follows is
explained
in detail in the book by Braides and Defranceschi \cite{BD} (see also
Dal Maso \cite {DM} and Buttazzo \cite{BU}).

Let $\Omega$ be a bounded subset of $\R^n$, let $p>1$ and let
$F_\e:\lp\orm\times\Ao\to [0,+\infty]$
be a family of functionals of the form
\begin{equation}\label{fel}
F_\e(u,U)=\cases{\displaystyle
\int_U f_{\e}(x,Du)\dx & if $u\in X_\e(U)$\cr\cr
+\infty & otherwise,
}
\end{equation}
for suitable function spaces $X_\e(U)$ and
$f_{\e}:\R^n\times\MM\to[0,+\infty)$
Borel functions. Suppose that there exist
 Borel functions $g_{\e}:\R^n\times\rr\to[0,+\infty)$,
convex and even in the second variable, with
\begin{eqnarray}\label{fgh}
g_\e(x,|F|)&\le& f_\e(x,F)\le C(1+g_\e(x,|F|))\le C(1+|F|^p),
\\
g_\e(x,2t)&\le& C(1+g_\e(x,t))\label{ggh}
\end{eqnarray}
for all $F\in\MM$, $x\in\Omega$ and $t\in\rr$. Growth conditions (\ref{fgh})
and (\ref{ggh}) are designed to include functions of the type $a_\e(x)|F|^p$ with
the only assumption $a_\e\ge 0$, thus allowing for zones where $a_\e=0$.
In the next section $a_\e$ will be the characteristic function of a 
set with fast-oscillating boundary.
Note that a general theory for functions satisfying
$$
0\le f_{\e}(x,F)\le C(1+|F|^p)
$$
only has not be developed yet. The aim of the direct method
of $\Gamma$-convergence is to prove a compactness result for the family
$(F_\e)$, giving a representation of the limit, and, possibly, complete the
description in terms of `homogenization formulas'.

{\it Step $1$}\ \  With fixed $(\e_j)$ extract a subsequence (not relabeled)
such that $F_\e(\cdot,U)$ $\Gamma$-converges to a functional $F_0(\cdot,U)$ for
all $U$ in a dense family of open sets $\cal U$ 
(see \cite{BD} Proposition 7.9);

{\it Step $2$}\ \ Thanks to (\ref{fgh}) and (\ref{ggh}), prove that
$F_0(u,\cdot)$ is the restriction of a
finite Borel measure to $\cal U$ for all $u\in\wup\orm$, so that by inner
regularity we
indeed have that $F_\e(\cdot,U)$ $\Gamma$-converges to a functional
$F_0(\cdot,U)$ on $\wup\orm$ for all $U\in \Ao$.
In this step is crucial the so-called  {\it fundamental
$\lp$-estimate}: for all $U, Y, Z \in {\cal A}(\Omega)$ with
$Y\subset\subset U$, and for all $\sigma>0$, there exists $M>0$ such
that for all $u,v\in \wup(\Omega;\R^m)$ one may find a cut-off
function
$\varphi \in C_0^{\infty}(U;[0,1])$, $\varphi= 1$ in $Y$, such that
\begin{eqnarray}\label{fund}\nonumber
F_\e(\varphi u+(1-\varphi)v, Y\cup
Z)&\le& (1+\sigma)(F_\e(u,U)+F_\e(v,Z))\\ && +
M\int_{(U\cap Z)\setminus Y}|u-v|^p\dx+\sigma\,.
\end{eqnarray}
Moreover, by again using the fundamental $\lp$-estimate it can be proven
that if $u\in \wup\orm \cap X_\e(U)$ for all $\e$ and
$F_0(u,U)<+\infty$ then there exist a sequence $u_{\e}\in X_\e(U)$ such that
$$
\lim_{\e\to 0} F_{\e}(u_{\e},U) = F_0(u,U)
$$
and $u_{\e}=u$ on a neighbourhood of $\partial U$
(see \cite{BD} Chapter 11);

{\it Step $3$}\ \ By the locality and semicontinuity properties of
$\Gamma$-limits
and by Step $2$ we can find a function $\varphi:\Omega\times\MM\to[0,+\infty)$
such that $0\le \varphi(x,F)\le C(1+|F|^p)$ and $F_0(u,U)=F_\varphi(u,U)$
for all
$u\in\wup\orm$
and $U\in\Ao$,
where
$$
F_\varphi(u)=\int_\Omega \varphi(x,Du)\dx.
$$
In the proof of this step a crucial point is the passage from the identity
$F_0(u)=F_\varphi(u)$ when $u$ is piecewise affine to a general $u$ by the
continuity of
$F_\varphi$ with respect to a convergence in which piecewise-affine functions
are dense
({\it e.g.} the strong $\wup$-convergence) (see \cite{BD} Chapter 9);

{\it Step $4$}\ \ If $f_{\e}(x,F)=f(\xe, F)$ with $f$ $1$-periodic in the first
variable then
by the periodicity of $f$ we deduce that $\varphi=\varphi(F)$
(see \cite{BD} Proposition 14.3);

{\it Step $5$}\ \ If $ g_{\e}(x,F)=g(\xe, F) $ with $g$ $1$-periodic in
the first variable then we
consider the auxiliary functionals
\begin{equation}\label{gel}
G_\e(u,U)=\cases{\displaystyle
\int_U g_\e(x,Du)\dx & if $u\in X_\e(\Omega)$\cr\cr
+\infty & otherwise.
}
\end{equation}
By Step $1$--$4$ we can assume that a function $\psi$ exists such that
$G_\e(\cdot,U)$ $\Gamma$-converges to the functional
$F_\psi(\cdot,U)$ on $\wup\orm$ for all $U\in\Ao$;

{\it Step $6$}\ \ Note that $\psi$ is convex. By an argument of
approximation by convolution prove that indeed the functional
$G_\e(\cdot,U)$ $\Gamma$-converges to the functional $F_\psi(\cdot,U)$ on
$\wuu\orm$ for all $U\in\Ao$. Define the `domain' of $F_\psi(\cdot,\Omega)$:
${\rm W}^{1,\psi}\orm=\{u\in\wuu\orm:\ F_\psi(u,\Omega)<+\infty\}$
(see \cite{BD} Theorem 14.8);

{\it Step $7$}\ \ Repeat Step $2$ and $3$ substituting the space $\wup\orm$
by the space
${\rm W}^{1,\psi}\orm$ thus obtaining the representation $F_0=F_\varphi$
on $\wuu\orm$;

{\it Step $8$}\ \ Deduce that $\varphi$ and $\psi$ do not depend on $(\e_j)$
by proving a homogenization formula (see \cite{BD} Proposition 21.12);

{\it Step $9$}\ \ Finally, the representation of $F_0$ on the whole
$\lp\orm$, and not only on $\wuu\orm$, can be
obtained in some cases by a more accurate study of the properties of $\varphi$.

\smallskip
We will have to modify Steps 1--9 above as to cover the case when
the domain of the limit is a `degenerate Sobolev Space'. In particular,
since the function $\psi$ obtained as in Step $5$ will be degenerate, a
suitable weighted Sobolev Space will have to be defined, which takes the
place of $\wuu\orm$ in Step $6$ above.
Moreover, we will have to deal with the fact that our
functions $f_{\e}, g_{\e}$ may be periodic only
in some variables, so that Step 8 will be harder to verify.
We will include all the details of the reasonings which do not fall directly
in this scheme, while we will feel free to refer to \cite{BD}
for those procedures which have become customary. 

It is worth mentioning that in some cases the arguments outlined above can be 
simplified by using some techniques (as blow-up arguments or the
theory of Young measures) that avoid to use the complex
localization procedure. As our problem is concerned those
methods seem harder to apply since the energies we consider are
coercive only on wildly oscillating sets.

\sect{Homogenization of media with oscillating profile}
Let $f: \R^{n-1}\mapsto [0,1]$ be a $1$-periodic
lower semicontinuous function and
$0\le \min f \le \sup f=1$, let $W:\R^{n-1}\times \MM \mapsto [0,+\infty)$
be a Borel function  $1$-periodic in the first variable satisfying
\begin{equation}\label{gc1}
\gamma |F|^p\le W(\x,F)\le \beta(1+|F|^p)
\end{equation}
for all $\x\in \R^{n-1}$ and $F\in \MM$, for some $1< p<+\infty$,
$0<\gamma\le\beta$. The set $\omega$ will be a fixed bounded open 
subset of $\R^{n-1}$ with Lipschitz boundary and $\Omega=\omega\times(-1,1)$.

%\vspace{-1cm}
\vspace{-0.5cm}
\begin{figure}[h]
\centerline{\psfig{figure=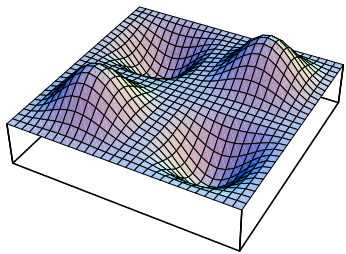,width=0.41\textwidth,angle=0}}
\vspace{-0.8cm}
\caption{the graph of a typical $f$ in the unit cell}
\end{figure}
\vspace{0.5cm}

In this section we compute the $\Gamma$-limit of functionals of the form
\begin{equation}\label{1}
J_\e(u)=\cases{\displaystyle
\int_{\Omega_\e} W\Bigl({x_\alpha\over\e},Du\Bigr)\dx
& if $u_{|\Omega_\e}\in\wup(\Omega_\e;\rr^m)$\cr\cr
+\infty & otherwise,}
\end{equation}
where
\begin{equation}\label{sd1}
\Omega_\e=\{x\in\Omega: \ |x_n|< f(x_\alpha/\e)\}.
\end{equation}

%\vspace{-1cm}
\vspace{-0.5cm}
\begin{figure}[h]
\centerline{\psfig{figure=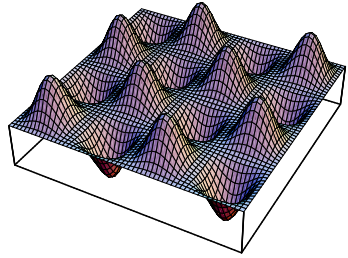,width=0.41\textwidth,angle=0}}  
\vspace{-0.8cm}
\caption{the upper profile of $\Omega_\e$ with $f$ as in Figure 1}
\end{figure}
\vspace{0.5cm}

The $\Gamma$-limit theorem will be stated and proved at the end of the
section after some preliminary results, which are needed to define
the domain of the $\Gamma$-limit and to explain the homogenization formula.

In orded to apply the method described in the previous section we 
introduce the localized version of the functionals $J_\e$: for all
$U$ open subset of $\Omega$ we define
\begin{equation}\label{1bbis}
J_\e(u,U)=\cases{\displaystyle
\int_{\Omega_\e\cap U} W\Bigl({x_\alpha\over\e},Du\Bigr)\dx
& if $u_{|\Omega_\e\cap U}\in\wup(\Omega_\e\cap U;\rr^m)$\cr\cr
+\infty & otherwise,}
\end{equation}
so that $J_\e(u)=J_\e(u,\Omega)$.

\bigskip
The first proposition contains the analog of Steps 1--4 of the direct method
of $\Gamma$-convergence as outlined in the previous section.

\begin{proposition}\label{gammacomp3}
From every sequence $(\e_j)$ of positive numbers converging to $0$
we can extract a subsequence (not relabeled) such that the $\Gamma$-limit
$$
J_0(u,U)=\Gamma\hbox{-}\lim_{j\to+\infty} J_{\e_j}(u,U)
$$
exists for all $u\in\wup(\Omega;\R^m)$ and $U$ open subsets of $\Omega$.
Moreover, there exists a Carath\'eodory function
$\varphi:(-1,1)\times\MM\to[0,+\infty)$
 such that
$$
J_0(u,U)=\int_U\varphi(x_n,Du)\dx
$$
for all $u\in\wup(\Omega;\R^m)$.
\end{proposition}

{\sc Proof.}
The functional $J_\e$ can be rewritten on
$X_\e(U)=\{ u\in \lp\orm \, :\, u_{|\Omega_{\e}\cap U}\in
\wup(\Omega_{\e}\cap U;\R^m)\}$ as
$$
J_\e(u,U)=\int_U \chi_{\Omega_\e}(x) W\Bigl({x_\alpha\over\e},Du\Bigr)\dx\,.
$$
We can then apply Steps 1-3 of Section 3 (see \cite{BD} Example 11.4 for the
proof of the $\lp$-fundamental estimate). Finally,
a translation argument in the $x_\alpha$-plane (completely analogous,
{\it e.g.}, to the one in the proof of \cite{BD} Proposition 14.3) shows that
$$
\int_{B_\rho(x_\alpha)\times(z-\eta,z+\eta)}\varphi(y,F)\dy
=
\int_{B_\rho(x'_\alpha)\times(z-\eta,z+\eta)}\varphi(y,F)\dy
$$
for all $\rho,\eta>0$, $x_\alpha,x'_\alpha$, $z$ such that
$$
\Bigl(B_\rho(x_\alpha)\times(z-\eta,z+\eta)\Bigr)
\cup
\Bigl(B_\rho(x'_\alpha)\times(z-\eta,z+\eta)\Bigr)
\subset \Omega\,.
$$
We then easily deduce  that $\varphi(x,F)=\varphi(x_n,F)$.
\qed

We will complete the proof of the homogenization theorem by characterizing 
the function $\varphi$ above (showing in particular 
that it does not depend on the sequence 
$(\e_{j})$), proving the existence of the $\Gamma$-limit $J_{0}$ on the whole 
$\lp\orm$ and showing that the integral representation in the previous 
proposition holds on the whole domain of $J_{0}$.
In order to get to this result,
we will have to define a number of auxiliary energies; 
here we streamline the organization of the rest of the 
section. First, in Section 4.1 we consider the case when $W(F)=\|F\|^p$.
We will denote by $\psi$ the function given by Proposition 
\ref{gammacomp3} corresponding to this particular choice of $W$.
For fixed $t$ the function $\psi(t,\cdot)$ is easily 
characterized by solving a $(n-1)$-dimensional (possibly, non coercive) 
homogenization problem.
It is possible then to define the `degenerate Sobolev Space' 
$\wup_{\psi}\orm$ of functions such that 
$\int_{\Omega}\psi(x_{n},Du)\dx<+\infty$, which turns out to
be the domain of the $\Gamma$-limit when $W(F)=\|F\|^p$, and hence
also in the general case by (\ref{gc1}).
In Section 4.2, in order to describe the function $\varphi$ 
in the general case, with fixed $t$ we consider the case when we replace 
the function $f$ by the characteristic function of $E_t=\{x_{\alpha}: 
f(x_{\alpha})>|t|\}$ (i.e., we deal with cylindrical domains). 
The function $\varphi(t,\cdot)$ will eventually
be given by the energy density of the corresponding $\Gamma$-limit.
Finally, in Section 4.3 we are able to consider general $W$ and $f$
and obtain the oscillating-boundary homogenization Theorem 4.15
as the consequence of the previous sections.

\subsection{An auxiliary problem. Definition of the limit domain}
In general, the limit functional $J_0$ exists and is finite also outside
$\wup(\Omega;\R^m)$. We first deal with the case of $J_0$ corresponding to
\begin{equation}\label{Fp}
W(x,F)=\|F\|^p,\hbox{ where } \|F\|^p=\sum_{j=1}^n|F_j|^p.
\end{equation}
By a careful description of the domain of the corresponding $\Gamma$-limit
we will identify the domain of $J_0$ as a suitable `degenerate Sobolev Space'
(see Definition \ref{root}) which, in view of the growth condition (\ref{gc1}),
will also be the domain of $J_0$ corresponding to energy densities other 
than (\ref{Fp}).

We recall a preliminary result.

\begin{theorem}\label{lessdim}
Let $E$ be a $1$-periodic set in $\R^N$\ie such that $\chi_E$ is a 
$1$-periodic function, and let
\begin{equation}
J^E_\e(v,U)=\cases{\displaystyle \int_{U\cap\e E } \|Dv\|^p\dx
& if $v_{|U\cap\e E }\in\wup(U\cap\e E ;\R^m)$\cr\cr
+\infty & otherwise.}
\end{equation}
Then the $\Gamma$-limit
$$
J^E_\Hom(v,U)=\Gamma\hbox{-}\lim_{\e\to 0} J^E_{\e}(v,U)
$$
exists for all $U$ bounded open subsets of $\R^N$ and $v\in\wup(U;\R^m)$.
Moreover, we have
$$
J^E_\Hom(v,U)=\int_U \varphi^E_\Hom(Dv)\dx
$$
for all $u\in\wup(U;\R^m)$, where $\varphi^E_\Hom$ is a positively homogeneous
function
of degree $p$, satisfying the formula
$$
\varphi^E_\Hom(F)=\inf\Bigl\{\int_{E\cap (0,1)^N}
\|Dv+F\|^p\dx: v
\in\wup_\loc(E;\R^m),\ 1\hbox{-periodic}\Bigr\}\,.
$$
\end{theorem}

{\sc proof.}
This theorem is a particular case of \cite{BD} Theorem 14.8,
the positive homogeneity of $\varphi^E_\Hom$ easily following from its
 definition.
\qed

For all $t\in(-1,1)$ we define
$$
\varphi_\#(t,\overline F)= \varphi^{E_t}_\Hom(\overline F)\,,
$$
the latter function being that given by the previous theorem,
with $N=n-1$ and $E=E_t=\{x_\alpha: \ f(x_\alpha)>|t|\}$.
We define also
\begin{equation}\label{psi}
\psi(t,F)= \varphi_\#(t,\overline F)+ {\cal
L}_{n-1}(E_t\cap(0,1)^{n-1})|F_n|^p.
\end{equation}

\begin{theorem}\label{2}
If $W=\|F\|^p$ and $\varphi$ is given by Proposition {\rm\ref{gammacomp3}}
then we have
$$
\varphi(t,F)=\psi(t,F)\,.
$$
In particular $\varphi$ does not depend on $(\e_j)$.
\end{theorem}

{\sc Proof.}
Let $(x,F)$ be such that $x_n$ is a Lebesgue point for $\varphi(\cdot, F)$.
Then
\begin{eqnarray}\label{A501}
\varphi(x_n,F)&=&
\lim_{\rho\to 0^{+}} \ave_{B_\rho(x_\alpha)\times(x_n-\rho,x_n)}
\varphi(y_n,F)\dy
\\
&=&
\nonumber
\lim_{\rho\to 0^+}
{J_0(Fy,B_\rho(x_\alpha)\times(x_n-\rho,x_n))
\over
|B_\rho(x_\alpha)\times(x_n-\rho,x_n)|}\,.
\end{eqnarray}
We consider the case  $x_n>0$ only, the case $x_n<0$ being dealt with 
using a symmetric
argument. Note that for $0<t<s<1$ we have $E_s\subseteq E_t$.
Let $u_j\to 0$ with
$u_j\in\wup_0(B_\rho(x_\alpha)\times(x_n-\rho,x_n)\cap \Omega_{\e_j})$
be such that
$$
J_0(Fy,B_\rho(x_\alpha)\times(x_n-\rho,x_n))
=\lim_{j\to+\infty} J_{\e_j} (Fy+u_j ,B_\rho(x_\alpha)\times(x_n-\rho,x_n))\,.
$$
Then,
\begin{eqnarray*}
&&J_{\e_j} (Fy+u_j ,B_\rho(x_\alpha)\times(x_n-\rho,x_n))\\
&=&
\int_{x_n-\rho}^{x_n}\int_{B_\rho(x_\alpha)}
\chi_{E_{y_n}}\Bigl({y_\alpha\over\e_j}\Bigr) 
\Vert \F+\D u_j\Vert^p dy_\alpha\, dy_n
\\
&&\quad+\int_{B_\rho(x_\alpha)}
\int_{x_n-\rho}^{x_n}\chi_{E_{y_n}} \Bigl({y_\alpha\over \e_j}\Bigr)
|F_n+D_n u_j|^p dy_n\,dy_\alpha
\\
&\ge &
\int_{x_n-\rho}^{x_n}\int_{B_\rho(x_\alpha)}
\chi_{E_{x_n}}\Bigl({y_\alpha\over\e_j}\Bigr) \Vert\F +\D u_j\Vert^p 
dy_\alpha\, dy_n
\\
&&\quad+\rho \int_{B_\rho(x_\alpha)}\chi_{E_{x_n}} 
\Bigl({y_\alpha\over \e_j}\Bigr) |F_n|^p \, dy_\alpha
\end{eqnarray*}
by Jensen's inequality. By using the lower limit inequality for the
$\Gamma$-convergence
in Theorem \ref{lessdim} with $E=E_{x_n}$, and by an
application of Fatou's Lemma, we get
$$
J_0(Fy,B_\rho(x_\alpha)\times(x_n-\rho,x_n))
\ge \rho\int_{B_\rho(x_\alpha)}\varphi_\#(x_n,\overline F)dy_\alpha
$$
$$
\qquad+\rho{\cal L}_{n-1}(B_\rho(x_\alpha))|F_n|^p{\cal L}_{n-1}(E_{x_n}
\cap(0,1)^{n-1})\,.
$$
Letting $\rho\to 0^+$ we obtain then by (\ref{A501})
$$
\varphi(x_n,F)\ge \varphi_\#(x_n,\overline F)+{\cal L}_{n-1}(E_{x_n}\cap(0,1)^
{n-1})|F_n|^p.
$$

Vice versa, let $v_j\to 0$ be such that $\overline Fy_\alpha+v_j(y_\alpha)$
is a recovery sequence for
$J^{E_{x_n}}_\Hom(\overline Fy_\alpha,
B_\rho(x_\alpha))$ along the sequence $(\e_j)$, and set
$$
u_j(y)=Fy+(v_j(y_\alpha),0)= (\overline Fy_\alpha+v_j(y_\alpha), F_ny_n)\,.
$$
We then have
\begin{eqnarray*}
&&\hskip-1cm
\int_{B_\rho(x_\alpha)\times(x_n, x_n+\rho)}\varphi(y_n,A)\dy\\
&\le &
\liminf_{j\to+\infty}J_{\e_j}(u_j,B_\rho(x_\alpha)\times(x_n, x_n+\rho))\\
&\le &
\liminf_{j\to+\infty}\int_{B_\rho(x_\alpha)\times(x_n, x_n+\rho)}
\chi_{E_{x_n}}\Bigl({y_\alpha\over\e_j}\Bigr)
\|Du_j\|^p\dy
\\
&= &\lim_{j\to+\infty}\rho\int_{B_\rho(x_\alpha)}\chi_{E_{x_n}}
\Bigl({y_\alpha\over\e_j}\Bigr) (\|\overline F+D_\alpha v_j\|^p
+|F_n|^p
)dy_\alpha
\\
&=&
\rho\int_{B_\rho(x_\alpha)} \varphi_\#(x_n,\overline F)dy_\alpha
\\&&\qquad+
\rho{\cal L}_{n-1}(B_\rho(x_\alpha))|F_n|^p{\cal L}_{n-1}(E_{x_n}
\cap(0,1)^{n-1}\,,
\end{eqnarray*}
which gives the missing inequality by (\ref{A501}).
\qed

\begin{remark}{\rm
With fixed $t$, we define the `kernel' of $\varphi_\#(t,\cdot)$ as
$$
\hbox{\rm Ker}\,\varphi_\#=\{\varphi_\#(t,\cdot)=0\}.
$$
Then $\hbox{\rm Ker}\, \varphi_\#$ is a linear space and its
dimension  is a multiple integer of $m$\ie
$$\hbox{\rm dim} \,\hbox{\rm Ker}\, \varphi_\# = km
\qquad \qquad  \hbox{\rm for some} \quad k=0,\ldots,n-1$$
and there exist  $\xi_{k+1},\ldots ,\xi_{n-1} \in \R^{n-1}$  such that
\[
\overline F =\left(\begin{array}{c}
F^1\\ \vdots\\F^m
\end{array}\right)
\in \hbox{\rm Ker}\, \varphi_\# \quad \Leftrightarrow \quad
\overline F \xi_i=0
\]
for each $i=k+1,\ldots,n-1$.
(Note that $k$ depends on $t$ fixed and $F^i$ denotes the
 $i$-th  row of
$\overline F$, $1\leq i \leq m$).

In fact,
since $\overline F \mapsto \varphi_\#(t,\overline F)$  is positively
homogeneous of degree $p$, convex and
even, $\hbox{\rm Ker}\, \varphi_\#$ is a linear space and satisfies
the following
properties:
if $\overline F \in  \hbox{\it Ker}\, \varphi_\#$ then

{\rm (i)} for each $(s_1,\ldots ,s_m) \in \Rm$
\[
 \left(\begin{array}{c}
s_1 F^1\\ \vdots \\s_m F^m
\end{array}\right) \in \hbox{\rm Ker}\, \varphi_\#;
\]

{\rm (ii)} $P\overline F \in \hbox{\rm Ker}\, \varphi_\#$
 for each permutation matrix $P\in \M^{m\times m} $.

Properties (i) and (ii) imply that if we  fix $F^1$ we can construct $m$
matrices linearly independent
\[
\left(\begin{array}{c}
F^1\\ 0\\ \vdots \\0
\end{array} \right),
\left( \begin{array}{c}
0\\ F^1 \\ \vdots \\ 0
\end{array}   \right),
\cdots,
\left( \begin{array}{c}
0\\ 0 \\ \vdots \\ F^1
\end{array}   \right)
\in \hbox{\rm Ker}\, \varphi_\#
\]
which span a subspace $\langle F^1 \rangle$  of $\hbox{\rm Ker}\, \varphi_\#$
of dimension $m$.

Now, if $\langle F^1 \rangle\neq {\rm Ker}\, \varphi_\#$,
we can single out a non-zero matrix in $\hbox{\rm Ker}\, \varphi_\#$
orthogonal to $\langle F^1 \rangle$, and, by using the same argument as above
taking its first row vector,
find other $m$ matrices which, together with the matrices constructed before,
form a linearly independent family.

By proceeding in this way, we end up with
$\eta_1,\ldots ,\eta_k \in \R^{n-1}$ such that for all $A\in \hbox{\rm Ker}\,
\varphi_\#$
 $$A^i=\sum_{j=1}^{k} s_{ij}\eta_j \qquad \qquad i=1,\ldots,m$$
with $s_{ij}\in \R$,
 which  means that the
$\hbox{\rm dim}\, \hbox{\rm Ker}\, \varphi_\#=km$ for some 
$k\in \{1,\ldots,n-1\}$.

The orthogonal subspace to $\langle \eta_1,\ldots ,\eta_k\rangle$ is
a vector subspace
of $\R^{m(n-1)}$ $\langle \xi_{k+1},\ldots ,\xi_{n-1}\rangle$
and the vectors of the two basis  satisfy, by definition,
 the conditions
$$\eta_i \xi_j=0 \qquad i=1,\ldots,k \qquad j=k+1,\ldots,n-1.$$
Hence, we can conclude that there exist vectors
$\xi_{k+1},\ldots ,\xi_{n-1}\in \R^{n-1}$ such that
$\overline F\in \hbox{\rm Ker}\, \varphi_\#$ if and only if
$\overline F \xi_i=0$
for each $i=k+1,\ldots,n-1$.

\bigskip

Since  $t\mapsto \varphi_\#(t,\overline F)$ is decreasing on $(0,1)$ and it is
coercive on $(0,\min f)$, there exist
$0\le \min f\le t_1\le\ldots\le t_k\le t_{k+1}\le\ldots\le t_{n-1}\le 1$
and $\xi_{k+1},\ldots ,\xi_{n-1} \in \R^{n-1}$   such that

{\rm(i)} $\varphi_\#(t,\overline F)$ is coercive on $(0,t_1)$;

{\rm(ii)} for each $k=1,\ldots,n-2$
$\varphi_\#(t,\overline F)= 0$ if and only if $\overline F\xi_i=0$ \linebreak
\indent for $i=k+1,\ldots,n-1$ on $(t_k,t_{k+1})$;

{\rm(iii)} $\varphi_\#(t,\overline F)= 0$ on $(t_{n-1},1)$.
}\end{remark}

\begin{definition}\label{root}
We define the `degenerate weighted Sobolev Space' $\wup_\psi(\Omega;\R^m)$
as the space of functions $u\in\lp(\Omega;\R^m)$ such that

{\rm (i)}  $D_n u\in\lp_\loc(\Omega;\R^m)$;

{\rm(ii)} $D_{(\xi _i,0)} u\in \lp_\loc(\omega\times(-t_i,t_i);\R^m)$
for $i=1,\ldots,n-1$;

{\rm (iii)} if \quad $\Phi:\Omega\to \M^{m\times (n-1)}$ \quad is
any measurable function such that \linebreak
\indent $\Phi \xi_i = D_{(\xi_i,0)}u \in \lp_\loc(\omega\times(-t_i,t_i);\Rm)$
for $i=1,\ldots,n-1$, then
$$
\int_\Omega \psi(x_n, \Phi | D_nu)\dx<+\infty\,.
$$

{\rm Clearly, the last integral is independent of the choice of $\Phi$; 
hence, it
will
be denoted by
$$
\int_\Omega \psi(x_n, Du)\dx\,,
$$
with a slight abuse of notation.}
\end{definition}

\begin{remark}{\rm
Note that in dimension $3$ ({\it i.e.}, $n=3$) the representation
of the space $\wup_\psi(\Omega;\R^m)$ is
particularly simple as, up to a rotation, we can assume that $\xi=e_2$.
In this case,
$\wup_\psi(\Omega;\R^m)$ is the space of functions
$u\in\lp(\Omega;\R^m)$ such that

{\rm(i)} $D_3 u\in\lp_\loc(\Omega;\R^m)$;

{\rm(ii)} $D_2 u\in \lp_\loc(\omega\times(-t_2,t_2);\R^m)$;

{\rm(iii)} $D_1 u\in \lp_\loc(\omega\times(-t_1,t_1);\R^m)$;

{\rm(iv)} if  $\Phi:\Omega\to {\Bbb M}^{m\times 2}$, $\Phi=(\Phi_1, \Phi_2)$ is
any measurable function such that \linebreak
\indent $\Phi_2 =D_2 u$ in $\omega\times(-t_2,t_2)$ and $\Phi_1 =
D_1 u$ in $\omega\times(-t_1,t_1)$, then
$$
\int_\Omega \psi(x_3, \Phi | D_3u)\dx<+\infty\,.
$$
}\end{remark}

\begin{example}\label{exxxon} {\rm
If $n=3$ and
$$
f(x_1,x_2)={1\over 2}+ {1\over 2}\sin^2(x_1)\sin^2(x_2)\,,
$$
then $\varphi_\#(t,\overline F)=\|\overline F\|^p$ if $|t|<1/2$ and $0$
otherwise, so that
$t_1=t_2=1/2$, and $\xi$ is any vector.
If instead
$$
f(x_1,x_2)={1\over 2}+ {1\over 2}\sin^2(x_1)\,,
$$
then $t_1=1/2$, $t_2=1$ and $\xi=(0,1)$.
}\end{example}

%\eject
%\vspace{-1cm}
\vspace{-0.5cm}
\begin{figure}[h]
\hbox{\hspace{0cm}\psfig{figure=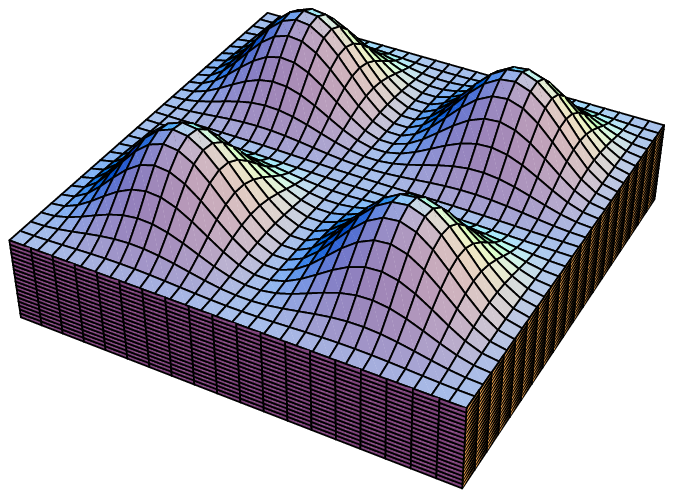,width=0.41\textwidth}\hspace{1cm}
\psfig{figure=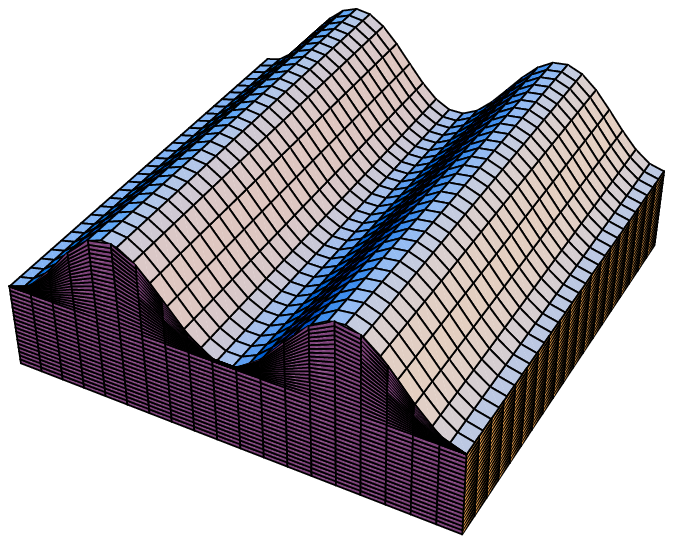,width=0.41\textwidth}}
\vspace{-0.8cm}
\caption{the oscillating profiles in Example 4.7}
\end{figure}
\vspace{0.5cm}

By using a convolution argument, we can improve Proposition \ref{gammacomp3} 
to give a characterization of the $\Gamma$-limit on the whole $\wup_\psi\orm$
and independent of the sequence $(\e_j)$. This result corresponds to
Step 6 in Section 3, and its proof uses the convexity of $F\mapsto \|F\|^p$
in an essential way.

\begin{proposition}\label{gammacomp4} Let $W=\|F\|^p$, and let $U$ be a
open subset of $\Omega$. Then 

{\rm(i)} if $u\in \lp(U;\R^m)\setminus \wup_\psi(U;\R^m)$ then there
exists the $\Gamma$-limit
$$
J_0(u,U)=\Gamma\hbox{-}\lim_{\e\to 0} J_{\e}(u,U)=+\infty;
$$

{\rm(ii)} if $u\in\wup_\psi(\Omega;\R^m)$ then there
exists the $\Gamma$-limit
$$
J_0(u,U)=\Gamma\hbox{-}\lim_{\e\to 0} J_{\e}(u,U)=
\int_U\psi(x_n,Du)\dx.
$$
\end{proposition}

{\sc Proof.}
We only outline the proof, as it closely follows that of \cite{BD} Theorem 
14.8, and details can be found therein.

Fix $u\in\lp(\Omega;\R^m)$ and $U$ an open subset
of $\Omega$. In order to compute $J_0(u,U)$ it is sufficient to
show that from every sequence $(\e_{j})$ we can extract a subsequence 
$(\e_{j_{k}})$ 
such that the $\Gamma$-limit along $(\e_{j_{k}})$ exists and
is independent of the subsequence.

We fix a sequence $(\e_{j})$.
By Theorem 4.3 the thesis of Proposition 4.1 holds
with $\psi$ in the place of $\varphi$. 
Upon possibly extracting a further 
subsequence, we may also assume that there exists the limit
$$
J_{0}(u,U)=\Gamma\hbox{-}\lim_{{j\to+\infty}}J_{\e_j}(u,U).
$$
Let $(\rho_j)$ be a sequence of mollifiers with
spt$\rho_j\subset B(0,{1\over j})\subset \R^{n-1}$, and define
$$
u_j(x)=\int_{B(0,{1\over j})} \rho_j(y)
u(x_{\alpha}-y,x_n)\dy.
$$

By the convexity of $J_{0}$ and its tran\-sla\-tion-in\-var\-iance
properties, we have 
$J_{0}(u_{j},U')\le J_{0}(u,U)$
for all $U'\subset\subset U$ such that
$U'\subset(y,0)+U$ for all $y\in\hbox{spt}\rho_j$.
By the convexity of $\psi$ the functional
$v\mapsto\int_{U'}\psi(x_n,Dv)\dx$ (if $v\in\lp(U';\rr^{m})\setminus\wup_{\psi}
(U';\rr^{m})$ this integral is set equal to $+\infty$) 
is lower semicontinuous with 
respect to the $\lp(U';\rr^{m})$ convergence. Hence,
we have
\begin{eqnarray*}
\int_{U'}\psi(x_{n},Du)\dx\le \liminf_{j\to+\infty}\int_{U'}\psi(x_{n},Du_{j})\dx
\le J_{0}(u,U).
\end{eqnarray*}
By the arbitrariness of $U'$ we get
\begin{equation}\label{abel}
\int_{U}\psi(x_{n},Du)\dx\le J_{0}(u,U),
\end{equation}
and in particular that $J_{0}(u,U)=+\infty$ if 
$u\in\lp(U;\rr^{m})\setminus\wup_{\psi}(U;\rr^{m})$,
so that (i) is proved.

Let now $u\in\wup_\psi\orm$. We first assume that 
$U\subset\subset U'\subset\subset\Omega$. 
By using the lower semicontinuity of 
$J_0$ and Jensen's inequality, we have
\begin{eqnarray*}
J_0(u,U)&\le &\liminf_{j\to+\infty}J_0(u_{j},U)
=\liminf_{j\to+\infty}\int_{U}\psi(x_{n},Du_{j})\dx
\\
&\le& \liminf_{j\to+\infty}\int_{U}\int_{B(0,{1\over j})}\rho_{j}(y)
\psi(x_{n},Du(x-(y,0)))\dx\dy
\\
&=& \liminf_{j\to+\infty}\int_{B(0,{1\over j})}\rho_{j}(y)
\int_{U+(y,0)}\psi(x_{n},Du)\dx\dy
\\
&\le& \liminf_{j\to+\infty}\int_{B(0,{1\over j})}\rho_{j}(y)\dy
\int_{U'}\psi(x_{n},Du)\dx= \int_{U'}\psi(x_{n},Du)\dx.
\end{eqnarray*}
By the arbitrariness of $U'$ we then get
\begin{equation}\label{porr}
 J_{0}(u,U)\le\int_{U}\psi(x_{n},Du)\dx,
\end{equation}
so that (ii) follows by taking (\ref{abel}) into account.

Finally, for arbitrary $U$, note that if $u\in\wup_\psi(\Omega;\R^m)$
then it can be approximated by a sequence $(v_j)$ of functions in
$\wup\orm$ such that $\int_{\Om}\psi(x_n,Dv_j)\dx$ are equi-bounded 
(we may use {\it e.g.} the argument in the proof
of \cite{EG} Section 4.2 Theorem 3); hence, by the lower 
semicontinuity of $J''=\Gamma$-$\limsup_j J_{\e_j}$, we have
$J''(u)<+\infty$. This fact implies (as in {\it e.g.} \cite{BD} Section 11.2) 
that $J''$ is inner-regular; i.e.,
$$
J''(u,U)=\sup\Bigl\{J''(u,V): V\subset\subset U\Bigr\}\,.
$$
Since (ii) holds with $V$ in the place of $U$ we easily get the thesis.
\qed

The following proposition clarifies the structure of $\wup_\psi$, and
implies that
the restrictions of functions $u\in\wup_\psi(\Omega;\rr^m)$ to relatively
compact subsets of $\omega\times (t_k, t_{k+1})$ are characterized as those
functions having directional derivatives $D_{k+1},\ldots, D_n$ $p$-summable.

\begin{proposition}\label{4}
Let  $k=1,\ldots,n-2$ and $s\in(t_k,t_{k+1})$.
There exist two positive constants
$\alpha_k(s)$ and $\beta_k$ such that

\begin{eqnarray}\label{5}
\alpha_k(s)\, \Bigl(\sum_{i=k+1}^{n-1} |\overline F\xi_i|^p + |F_n|^p\Bigr)
\le \psi(t,F)\le \beta_k\,\Bigl(\sum_{i=k+1}^{n-1}  |\overline F\xi_i|^p
+ |F_n|^p\Bigr)
\end{eqnarray}
for all $F\in \MM$ and $t\in(t_k,s]$ .
\end{proposition}

{\sc Proof.}
Since $\overline F \mapsto \varphi_\#(t,\overline F)$  is positively
homogeneous of degree $p$ and convex, if  $t\in(t_k,t_{k+1})$
we easily deduce that
$$
\varphi_\#(t,\overline F)\le  c \sum_{i=k+1}^{n-1}
\varphi_\#(t,\Xi_i) |\overline F\xi_i|^p
$$
where
\[
\Xi_i=\left( \begin{array}{c}
\xi_i\\ 0 \\ \vdots \\ 0
\end{array}   \right).
\]

If we denote
$$
\beta_k'=\max_{i=k+1,\ldots,n-1}\sup_{t\in [0,1)}c\,
\varphi_\# (t,\Xi_i)
$$
then
\begin{eqnarray}\label{10}
\varphi_\#(t,\overline F)\le \beta_k'
\sum_{i=k+1}^{n-1} |\overline F\xi_i|^p.
\end{eqnarray}
On the other hand we have that

\begin{eqnarray*}
{\varphi_\#(t,\overline F)\over  \sum_{i=k+1}^{n-1}  |\overline F\xi_i|^p}
&\ge&
c\, {\varphi_\#(t,(\overline F\xi_{k+1},\ldots ,\overline F \xi_{n-1}))\over
\|(\overline F\xi_{k+1},\ldots ,\overline F\xi_{n-1})\|^p}
\\
&\ge&
c \inf\{
\varphi_\#(t,G)\, : \,  G\in S^{n-1}\cap {\hbox{\rm Ker}\varphi_\#}^\bot\}
\end{eqnarray*}
by p-homogeneity.
Note that
$t \mapsto c\, \inf\{
\varphi_\#(t,G)\, : \, G\in S^{n-1}\cap {\hbox{\rm Ker}\varphi_\#}^\bot\}=
c(t) $
is decreasing on $(0,1)$  and
$$\inf_{t\in (t_k,s]} c(t)=\alpha'_k(s)>0,$$
so that we get
\begin{eqnarray}\label{11}
\varphi_\#(t,\overline F) \ge
\alpha'_k(s)\, \sum_{i=k+1}^{n-1} |\overline F\xi_i|^p.
\end{eqnarray}
Let
$$
\alpha_k(s)=\min\{\alpha'_k(s),\inf_{t\in (t_k,s]}{\cal L}_{n-1}
(E_t\cap(0,1)^{n-1})\}
$$
and
$$
\beta_k=\max\{\beta'_k,1\},
$$
then (\ref{5}) follows by Theorem \ref{2}, (\ref{10}) and (\ref{11}).
\qed

\begin{proposition}\label{prprpr}
Fix $t\in (t_k,t_{k+1})$, for $k=0,\ldots,n-1$ ($t_0=0,t_n=1$).
If $\psi$ is given by {\rm(\ref{psi})} then
\begin{eqnarray*}
\psi(t,F)&=&\min \Bigl\{\int_{(0,1)^n\cap (E_t\times (0,1))}
\|Dw\|^p\dx :
\\ \nonumber&&\qquad\qquad
w\in\wup_\loc(E_t\times (0,1);\R^m),\
 w-Fx\ 1\hbox{-periodic}\Bigr\}\,.
\end{eqnarray*}

\end{proposition}

{\sc Proof.}
Let $w$ be a test function for the minimum problem above, then
\begin{eqnarray*}
&&\int_{(0,1)^n\cap (E_t\times (0,1))}
\|Dw\|^p\dx\\
&=&
\int_{(0,1)^n\cap (E_t\times (0,1))}
\|D_{\alpha}w\|^p\dx + \int_{(0,1)^n\cap (E_t\times (0,1))}
|D_nw|^p\dx \\
&\ge &
\int_0^1\min \Bigl\{\int_{E_t\cap (0,1)^{n-1}}
\|Dv\|^p\dx_{\alpha}:
\\ \nonumber&&\qquad\qquad\qquad
v\in\wup_\loc(E_t;\R^m),\ v-\overline{F}x_{\alpha}\
1\hbox{-periodic}\Bigr\}\dx_n
\\
&&
+\int_{E_t\cap (0,1)^{n-1}}\Bigl(\int_0^1 |D_nw|^p\dx_n  \Bigr)\dx_{\alpha}\\
&\ge &
\varphi_\#(t,\overline{F})+{\cal L}_{n-1}(E_t\cap(0,1)^{n-1}) |F_n|^p=
\psi(t,F)
\end{eqnarray*}
by Jensen's inequality and the description of $\varphi_\#$ 
(see Theorem \ref{lessdim}); hence, 
\begin{eqnarray*}
\psi(t,F)&\le &\min\Bigl\{\int_{(0,1)^n\cap E_t\times (0,1)}
\|Dw\|^p\dx:
\\ \nonumber&&\qquad\qquad
w\in\wup_\loc(E_t\times(0,1);\R^m),\ w-Fx\ 1\hbox{-periodic}\Bigr\}
\end{eqnarray*}
by Theorem \ref{2}.

Conversely, given a function $v$ such that $v-\overline{F}x_{\alpha}$ is 
$1$-periodic, we can construct a test function $w$, such that
$w-Fx\ 1\hbox{-periodic}$, as
$$w(x)=v(x_{\alpha})-F_nx_n\,.$$
We then have
\begin{eqnarray*}
&&\int_{(0,1)^n\cap (E_t\times (0,1))}
\|Dw\|^p\dx\\
&=&
\int_{(0,1)^n\cap (E_t\times (0,1))}
(\|\D v\|^p+|F_n|^p)\dx\\
&=&
\int_{E_t\cap (0,1)^{n-1}}
\|\D v\|^p\dx_{\alpha}+
{\cal L}_{n-1}(E_t\cap(0,1)^{n-1}) |F_n|^p\\
&\ge &
\min \Bigl\{\int_{(0,1)^n\cap (E_t\times (0,1))}
\|Dw\|^p\dx:
\\ \nonumber&&\qquad\qquad
w\in\wup_\loc(E_t\times (0,1);\R^m),\ w-Fx\ 1\hbox{-periodic}\Bigr\}
\end{eqnarray*}
and hence the converse inequality
\begin{eqnarray*}
\psi(t,F)&=&\min\Bigl\{\int_{E_t\cap (0,1)^{n-1}}
\|Dv\|^p\dx_{\alpha}:
\\ \nonumber&&\qquad\qquad
v\in\wup_\loc(E_t;\R^m),\ v-\overline{F}x_{\alpha}\ 1\hbox{-periodic}\Bigr\}\\
&&
+{\cal L}_{n-1}(E_t\cap(0,1)^{n-1}) |F_n|^p\\
&\ge &
\min \Bigl\{\int_{(0,1)^n\cap (E_t\times (0,1))}
\|Dw\|^p\dx:
\\ \nonumber&&\qquad\qquad
w\in\wup_\loc(E_t\times (0,1);\R^m),\ w-Fx\ 1\hbox{-periodic}\Bigr\}\,
\end{eqnarray*}
is obtained as desired.
\qed

Now we can turn our attention to the case with a general $W$.
Now that a natural domain for the limit functional is defined, we
can easily state and prove a
compactness result that partly improves Proposition \ref{gammacomp3}.

\begin{theorem}\label{8a}
Let $J_\e$ be given by {\rm(\ref{1bbis})}. Then
for every sequence $(\e_j)$ of positive numbers converging to $0$
there exists a subsequence (not relabeled) such that
the $\Gamma$-limit
$$
J_0(u,U)=\Gamma\hbox{-}\lim_{j\to+\infty} J_{\e_j}(u,U)
$$
exists for all $u\in\wup_\psi(\Omega;\Rm)$ and $U$ open subsets of $\Omega$.
Moreover $J_0(u,\cdot)$ is the restriction of a Borel measure to $\Ao$.
\end{theorem}

{\sc Proof.}
By (\ref{gc1}) and Proposition \ref{gammacomp4}
we deduce the condition
\begin{eqnarray}\label{8}
\Gamma\hbox{-}\limsup_{\e\to 0}J_{\e}(u,U)\le \beta \int_U(1+\psi(x_n,Du))\dx
\end{eqnarray}
if $u\in\wup_\psi(U;\Rm)$ and $U$ is an open subset of $\Omega$. Then, we can
follow the Steps 1--3 in Section 3 to prove the compactness of $(J_\e)$ and
that $J_0(u,\cdot)$ is the restriction  of a Borel measure to $\Ao$.
\qed

\subsection{Homogenization of cylindrical domains}
It remains now to extend the integral representation of
Proposition \ref{gammacomp3} and characterize its integrand.  We first
deal with the case of `cylindrical' domains\ie we consider $\chi_E$ in place 
of $f$, with $E$ a $1$-periodic open subset of $\R^{n-1}$.

Let $t_1,\ldots,t_{n-1}$ be the points in $(0,1)$ introduced to 
characterize the 
`degenerate weighted Sobolev Space' in Definition \ref{root}. Since in the
following we will choose $E=E_t$ ($E_t$ defined as 
$\{x_\alpha: \ f(x_\alpha)>|t|\}$) we introduce the following notation:
with fixed $t\in(0,1)$, $t\not= t_k$ for $k=1,\ldots,n-1$, consider
the set $E_t$ and the functional
\begin{equation}\label{7}
J_\e^t(u,U)=\cases{\displaystyle
\int_{\Omega_{\e}\cap U_{\e}}
 W\Bigl({x_\alpha\over\e},Du\Bigr)\dx
& if $u\in\wup(\Omega_{\e}\cap U_{\e};\R^m)$\cr\cr
+\infty & otherwise,}
\end{equation}
where $U_{\e}=U\cap (\e E_t\times (-1,1))$.
Note that the integrand of $J^t_\e$ satisfies the following growth conditions
\begin{eqnarray}\label{6}
\gamma g(x,A)\le \chi_{E_t\times (-1,1)}
W(x_\alpha,A) \le \beta (1+g(x,A))
\end{eqnarray}
where
$g(x,A)=\chi_{E_t\times (-1,1)}(x) \|A\|^p$
is obviously 1-periodic in $x$, convex in $A$ and satisfying
$$
0\le g(x,A) \le 1+\|A\|^p
\hbox{ and }
g(x,2A)\le c(1+g(x,A))$$
for all $A\in \MM$.

\begin{remark}{\rm Note that if we fix $t\in (t_{k-1},t_k)$ and consider
$\chi_{E_t}$ in place of $f$ then  $\wup_{\psi}(\Omega;\R^m)$ turns
out to be the space
\begin{eqnarray*}
\wup_k(\Omega;\R^m)&= &\{u\in \lp(\Omega;\Rm)\,
: D_n u \in \lp(\Omega;\R^m),
D_{\xi_i}u\in \lp(\Omega;\R^m)
\\
&&
\hskip 4cm i=k,\ldots,n-1\}
\end{eqnarray*}
if $k=1,\ldots,n-1$, and
$$
\wup_n(\Omega;\R^m)=
\{u\in \lp(\Omega;\Rm)\, :\,  D_n u \in \lp(\Omega;\R^m)\}.
$$
if $k=n$.
}\end{remark}

%\vspace{-1cm}
\vspace{-0.5cm}
\begin{figure}[h]
\hbox{\hspace{0cm}\psfig{figure=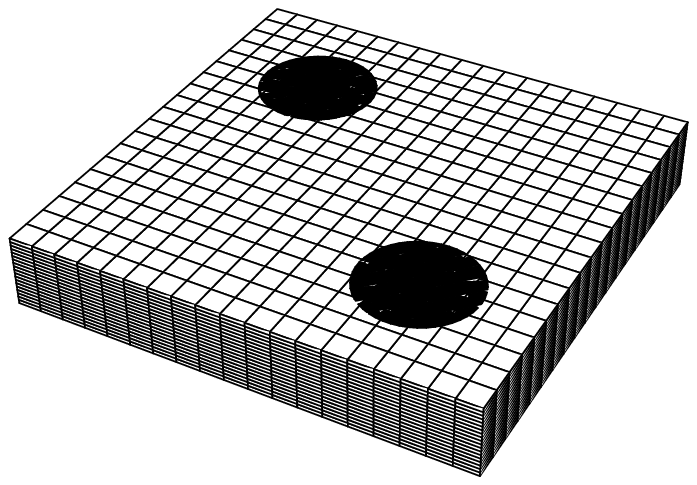,width=0.41\textwidth}\hspace{1cm}
\psfig{figure=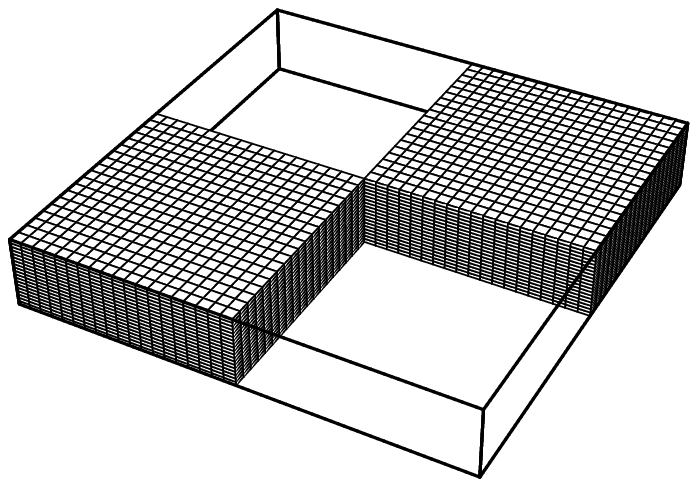,width=0.41\textwidth}}
\vspace{-1cm}
\centerline{\psfig{figure=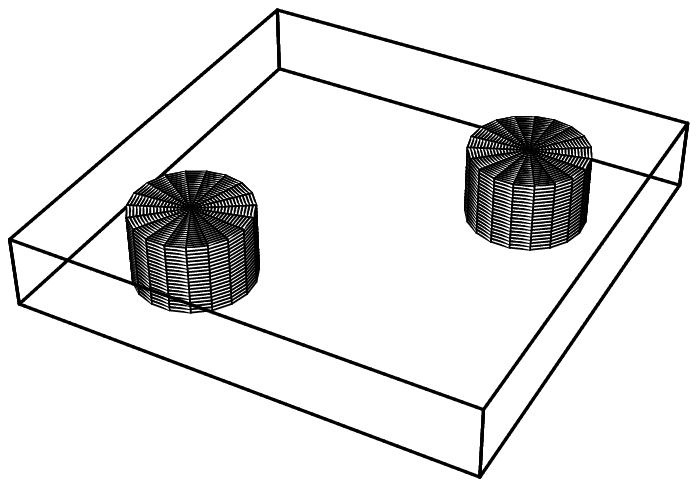,width=0.41\textwidth}}
\vspace{-0.8cm}
\caption{cylindrical domains $E_t \times (-1,1)$ related to the function $f$
 in Figure 1 for different values of $t$}
\end{figure}
\vspace{0.5cm}

\begin{theorem}\label{9}
Let $t\in (t_{k-1},t_k)$ and let $J_\e^t(\cdot,U)$ be defined by
{\rm(\ref{7})}. Then the $\Gamma$-limit
$$
J_0^t(u,U)=\int_U W^t_{\Hom}(Du)\dx
$$
exists for each $u\in \wup_k(\Omega;\Rm)$ and $U$ open subset of 
$\Omega$,
where $W^t_{\Hom}$ is given
by
\begin{eqnarray*}
W^t_{\Hom}(A)&=&
\lim_{T\to+\infty} \inf\Bigl\{{1\over T^n}\int_{(0,T)^n}
 \chi_{E_t}(x_\alpha) W(x_{\alpha},A+Du(x))\dx :
\\ \nonumber&&\qquad\qquad
u\in \wup_0((0,T)^n;\Rm)\Bigr\}
\end{eqnarray*}
for all $A\in \MM$.
\end{theorem}

{\sc Proof.}
By taking Theorem \ref{8a} into account with $\chi_{E_t}$ in the place 
of $f$,
and repeating word for word the proof of the integral representation
theorem \cite{BD} Theorem 9.1, replacing $\wup\orm$ by $\wup_k\orm$,
we obtain an integral representation on the whole $\wup_k\orm$.
The integrand of this representation must coincide with
the function $\varphi=\varphi(x_n,F)$ provided by Proposition
\ref{gammacomp3} with
$\chi_{E_t}$ in the place 
of $f$. Since the functionals are clearly invariant by
translations in the direction $x_n$ we have indeed $\varphi=\varphi(F)$.
To prove the asymptotic formula we can repeat the proof
of Proposition 21.12 in \cite{BD}.
\qed

\subsection{The general case}
We can eventually proceed to dealing with the general case.

\begin{proposition}\label{intrepgen} Let
$J_\e$ be given by {\rm(\ref{1bbis})}. Then
the $\Gamma$-limit
$$
J_0(u,U)=\Gamma\hbox{-}\lim_{\e\to 0} J_{\e}(u,U)
$$
exists for all $u\in\wup_\psi(\Omega;\R^m)$ and $U$ open subsets of $\Omega$.
Moreover, for such $u$ we have
$$
J_0(u,U)=\int_U \varphi(x_n,Du)\dx,
$$
where $\varphi$ is given by Proposition
{\rm\ref{gammacomp3}}.
\end{proposition}

{\sc Proof.}
We have to extend the representation of $J_0$ given by Proposition
\ref{gammacomp3} to $\wup_\psi\orm$.
Note that $\varphi$ is a Carath\'eodory function (see \cite{BD} Theorem
9.1, Step 3).
As explained in Step $3$ of Section $3$, a crucial argument used to obtain
an integral
representation result is the continuity in $\wup_\psi\orm$ of the functional
$$
u\mapsto \int_U \varphi (x_n,Du)\dx
$$
along some
strongly converging sequences of piecewise-affine functions. We only prove
this property, as the rest of the proof follows exactly that of
\cite{BD} Theorem 9.1 (Steps 1--3, 5 and 6; the proof below replaces Step 4).

\smallskip
Let $U=\bigcup_{k=0}^{n-1} U_k$ where $U_k\subset\subset
\omega\times(t_k,t_{k+1})$, ($t_0=0, t_n=1$);
we can find functions  $u_j\in \wup_\psi(\Omega;\Rm)$ such
that their restrictions to $U$ are piecewise affine and
$u_j,D_nu_j$ converge strongly to $u,D_nu$ in $\lp(U;\Rm)$, respectively,
while $D_{(\xi_i,0)}u_j$ converge strongly to $D_{(\xi_i,0)}u$ in 
$\lp(U_i;\Rm)$.

We will use some estimates deriving from the inequality
$\varphi(t,F)\le \beta(1+\psi(t,F))$,
which follows trivially from (\ref{gc1}).
By Proposition \ref{4} we have that
\begin{eqnarray*}
\psi(x_n,Du)&\le&\beta_k\Bigl(\sum_{i=k+1}^{n-1} |D_{(\xi_i,0)}u|^p +|D_nu|^p
\Bigr)
\\
\psi(x_n,Du_j)&\le&\beta_k\Bigl(\sum_{i=k+1}^{n-1} |D_{(\xi_i,0)}u_j|^p +
|D_nu_j|^p\Bigr)
\end{eqnarray*}
on $\omega\times (t_k,t_{k+1})$. Note that by (\ref{8})
\begin{eqnarray*}
 \int_U\varphi(x_n,Du_j)\dx &\le& \sum_{k=0}^{n-2} \int_{U_k\cap \omega\times
(t_k,t_{k+1})}
\beta \Bigl(1+\sum_{i=k+1}^N \beta_k |D_{(\xi_i,0)}u_j|^p\Bigr)\dx
\\
&&+
\beta\int_{U}\beta_k |D_nu_j|^p\dx.
\end{eqnarray*}
If we use the continuity of $\varphi$ in the second variable and
apply Fatou's lemma to the sequences
\begin{eqnarray*}
\beta\int_{U}\beta_k |D_nu_j|^p\dx
 &+& \sum_{k=0}^{n-2} \int_{U_k\cap \omega\times
(t_k,t_{k+1})}
\beta \Bigl(1+\sum_{i=k+1}^{n-1} \beta_k |D_{(\xi_i,0)}u_j|^p\Bigr)\dx
\\
&\pm&\int_U\varphi(x_n,Du_j)\dx\end{eqnarray*}
we get that
$$
\int_U\varphi(x_n,Du)\dx=\lim_{j\to +\infty} \int_U\varphi(x,Du_j)\dx.
$$
Hence, we have proved the integral representation for sets of the type
$U=\bigcup_{k=0}^{n-1} U_k$ where
$U_k\subset\subset \omega\times(t_k,t_{k+1})$.
A symmetric argument applies to the case where 
$U=\bigcup_{k=0}^{n-1} U_k$, with 
$U_k\subset\subset \omega\times(-t_{k+1},-t_k)$.
Since $J_0(u,\cdot)$ is a measure absolutely continuous with respect to
Lebesgue measure, we conclude that the integral representation holds for all open
subsets $U$ of $\Omega$.
\qed

Finally, the oscillating-boundary homogenization theorem reads as follows.

\begin{theorem}\label{oscihom}
Let $J_\e$ be given by {\rm(\ref{1})}. Then
the $\Gamma$-limit
$$
J_0(u)=\Gamma\hbox{-}\lim_{\e\to 0} J_{\e}(u)
$$
exists for all $u\in\lp(\Omega;\R^m)$, and we have
$$
J_0(u)=\cases{\displaystyle\int_\Omega W_{\Hom}(|x_n|,Du)\dx
 & if $u\in\wup_\psi(\Omega;\R^m)$\cr\cr
+\infty & otherwise,}
$$
where $W_{\Hom}(t,A)=W^t_{\Hom}(A)$  for a.e. $t\in (0,1)$,
and $W^t_{\Hom}$ is given by Theorem {\rm\ref{9}}.
Moreover, if $u\in\wup\orm$ there exists a family
$(u_\e)$ converging to $u$ in $\lp\orm$, such that $u-u_\e$ has compact 
support in $\Omega$ and
$J_0(u)=\lim_{\e\to 0}J_\e(u_\e)$.
\end{theorem}

{\sc Proof.}
It is sufficient to
compute the $\Gamma$-limit for $u\in\wup_\psi(\Omega;\R^m)$,
since by comparison with Proposition \ref{gammacomp4}(i)
we immediately have $J_0(u)=+\infty$ if $u\not\in\wup_\psi(\Omega;\R^m)$.
Let $\varphi$ be given by Proposition \ref{gammacomp3}; it
remains to prove that $\varphi$ satisfies an asymptotic formula.

Let $x_{n}>0$, let $0<\rho<x_n$ and consider the functionals (\ref{7}) with
$t=x_n-\rho$ and $t=x_n$ so that
\begin{eqnarray*}
&&\hskip-2cm
J_{\e}^{x_n-\rho}(Ax,(0,1)^{n-1}\times (x_n-\rho,x_n))\\
&\ge&
\int_{(0,1)^{n-1}\times (x_n-\rho,x_n)}
 \chi_{E_{y_n}}\Bigl({y_{\alpha}\over \e}\Bigr)
W\Bigl({x_{\alpha}\over \e},A\Bigr)\dy
\\
&\ge&
J_{\e}^{x_n}(Ax,(0,1)^{n-1}\times (x_n-\rho,x_n))\,.
\end{eqnarray*}
By Theorem \ref{9}
\begin{eqnarray*}
\rho\, W^{x_n-\rho}_{\Hom}(A)&\ge&
\Gamma\hbox{-}\lim_{\e\to 0} J_{\e}(Ax,(0,1)^{n-1}\times (x_n-\rho,x_n))
\\
&\ge&
\rho\, W^{x_n}_{\Hom}(A).
\end{eqnarray*}
Taking into account that
$$
\Gamma\hbox{-}\lim_{\e\to 0} J_{\e}(Ax,(0,1)^{n-1}\times (x_n-\rho,x_n))=
\int_{(0,1)^{n-1}\times (x_n-\rho,x_n)} \varphi(y_n,A)\dy
$$
we get
$$
W^{x_n-\rho}_{\Hom}(A)\le {1\over \rho} \int_{(x_n-\rho, x_n)}
\varphi(y_n,A)\dy_n \le
W^{x_n}_{\Hom}(A).
$$
Since $t\mapsto W^t_{\Hom}(A)$ and $t\mapsto \varphi(t,A)$ are decreasing
functions on $(0,1)$,
there exists a subset $M$ of $(0,1)$, $|M|=0$, such that they are
continuous
on $(0,1)\setminus M$; hence, by passing to the limit as $\rho \to 0$ we get
$$
\varphi(x_n,A)=W^{x_n}_{\Hom}(A)
$$
for every $x_n\in (0,1)\setminus M$. For $x_{n}<0$ it suffices to 
apply a symmetric argument.

The last statement follows by a well-known argument of stability of
$\Gamma$-con\-ver\-gence by compatible boundary data due to De Giorgi
(see \cite{BD} Section 11.3).
\qed

\sect{Thin films with fast-oscillating profile}
In this section we establish the second goal of the paper; that is,
to prove that the $\Gamma$-limit of functionals $E_{\e,\delta}$ as in 
(\ref{1.01})
when $\e\to 0$ and $\delta<<\e$, is given by first applying the theory
constructed in the previous section with $\e$ as a parameter and letting
$\delta\to 0$, and subsequently letting $\e\to 0$. The final result can
be summarized as follows, in a $n$-dimensional setting.

\begin{theorem}\label{tfft}
Let $f: \R^{n-1}\to [0,1]$ be a $1$-periodic lower semicontinuous
function with $0< \min f \le \sup f=1$, let $W:\MM \to [0,+\infty)$ be a
convex function satisfying
$$
\gamma |F|^p\le W(F)\le \beta(1+|F|^p)
$$
for all $F\in \MM$ and for some $1<p<+\infty$, $0<\gamma\le\beta$.
Let $\delta:(0,+\infty)\to(0,+\infty)$ be such that
$$
\lim_{\e\to 0}{\delta(\e)\over\e}=0.
$$
Let $\omega$ be a bounded open subset of $\rr^{n-1}$ and
let $\Omega_\e\subset \omega\times(-1,1)$ be defined by
\begin{equation}
\Omega_\e=
\Bigl\{
x\in\R^n:\ |x_n|< f\Bigl({x_\alpha\over\delta(\e)}\Bigr), x_\alpha\in\omega
\Bigr\}.
\end{equation} Define $E_\e:\lp(\omega\times(-1,1))\to [0,+\infty]$
by
\begin{equation}
E_{\e}(u)=\cases{\displaystyle
\int_{\Omega_\e} W\Bigl(D_\alpha u,{1\over\e}D_n u\Bigr)\dx
& if $u_{\bigl|\Omega_\e}\in\wup(\Omega_\e;\rr^m)$\cr\cr
+\infty & otherwise.}
\end{equation}
Then the $\Gamma$-limit as $\e\to0$ of $E_\e$ is given by 
\begin{equation}
E(u)=\cases{\displaystyle
\int_{\omega\times(-1,1)} \overline W_\Hom
(D_\alpha u)\dx & if $u\in\wup(\omega\times(-1,1);\rr^m)$ and $D_n u=0$
\cr\cr
+\infty & otherwise, }
\end{equation}
where $\overline W_\Hom:\M^{m\times (n-1)}\to[0,+\infty)$ is given by
\begin{equation}\label{treat}
\overline W_{\Hom}(\overline F)=
\int_0^1 \inf_{F_n}W_{\Hom}(t,\overline F|F_n)\dt,
\end{equation}
and $W_\Hom$ by
\begin{eqnarray}\label{true}\nonumber
W_{\Hom}(t,F)&=&
\inf\Bigl\{\int_{(0,1)^n}
 \chi_{E_t}(x_\alpha)\, W(F+Du(x))\dx :
\\ &&\qquad\qquad
u\in \wup_\loc(\rr^n;\Rm)\hbox{ $1$-periodic}\Bigr\}
\end{eqnarray}
for all $t\in(0,1)$ and $F\in \MM$, where $E_t=\{f>t\}$.
\end{theorem}

\subsection{Proof of the result}
In order to simplify the proof without loosing sight of the main intricacies
of the argument, 
we deal only with the case where $\e=1/j$ and $\delta=\e^2$. 
The general case
can be dealt with similarly, by introducing some error terms.
We define, with a slight abuse of notation,
$$
\Omega_k=\{x\in\Omega: \ |x_n|< f(k x_{\alpha})\}
$$
and for $k = j^2$, $j\in \NN$
$$
E_j(u,U)=\int_{\Omega_{j^2}\cap U}W(D_{\alpha}u|jD_nu)\dx
$$
for all $u_{|\Omega_{j^2}\cap U}\in\wup(\Omega_{j^2}\cap U;\rr^m)$.

By the compactness result Theorem 2.5 in \cite{BFF} we can suppose
that there exists $W_0: \M^{m\times (n-1)} \to [0, +\infty)$ such that 
$E_j(u,U) $ $\Gamma$-converge for all sets
of the form $U=U'\times(-1,1)$ or $U=U'\times(0,1)$ to the functional given by
\begin{equation}\label{inf}
E_0(u,U)=
\cases{\displaystyle
\int_{U}W_0(D_{\alpha}u)\dx
& if $u\in\wup(U;\rr^m)$ and  $D_n u=0$\cr\cr
+\infty & otherwise.}
\end{equation}

\begin{proposition}\label{gg}
For all $\overline F\in \M^{m\times (n-1)}$ define
\begin{eqnarray}\label{a}
\overline W_{\Hom}(\overline F)&=&\inf \Bigl\{\int_{(0,1)^n}
W_{\Hom}(x_n, Du+\overline F)\dx:
\\ \nonumber&&\qquad\qquad
u\in\wup_\loc(\R^n;\R^m),\
u\ 1\hbox{-periodic\, in}\,\, x_{\alpha}\Bigr\}\,.
\end{eqnarray}
Then
\begin{equation}\label{ccc}
\overline W_{\Hom}(\overline F)=
\int_0^1 \widetilde{W}_{\Hom}(t,\overline F)\dt,
\end{equation}
where
\begin{equation}\label{c}
 \widetilde{W}_{\Hom}(t,\overline F)=
\inf_{F_n}W_{\Hom}(t,\overline F|F_n)
\end{equation}
and $\overline F \mapsto  \widetilde{W}_{\Hom}(t,\overline F)$ is convex.
\end{proposition}

{\sc Proof.}
It can be easily proved that $\overline F \mapsto  \widetilde{W}_{\Hom}
(t,\overline F)$ is convex.

With fixed $\eta>0$, by the Measurable Selection Criterion
(see e.g \cite{CV}),
we can find $G_n(t)$ a measurable function such that
$$
W_{\Hom}(t,\overline F|G_n)\le \inf_{F_n}W_{\Hom}(t,\overline F|F_n)+
\eta.
$$
We can consider
$$
u(x_{\alpha},x_n)=\int_0^{x_n} G_n(s)\,ds
$$
as test function  in (\ref{a}). We then get
$$
\overline W_{\Hom}(\overline F)\le \int_0^1
W_{\Hom}(x_n,\overline F|G_n(x_n)) \dx_n
$$
and so
$$
\overline W_{\Hom}(\overline F)\le \int_0^1
\inf_{F_n}W_{\Hom}(t,\overline F|F_n) \dt +\eta=\int_0^1
 \widetilde{W}_{\Hom}(t,\overline F) \dt +\eta.
$$
Conversely,
\begin{eqnarray*}
\overline W_{\Hom}(\overline F)&\ge&
\inf \Bigl\{\int_{(0,1)^n}
\widetilde{W}_{\Hom}(x_n, D_{\alpha}u+\overline F)\dx:
\\ \nonumber&&\qquad
u\in\wup_\loc((0,1)^n;\R^m),\
u\ 1\hbox{-periodic\, in}\, x_{\alpha}\Bigr\}\\
&\ge&
\int_0^1\Bigl(
\inf \Bigl\{\int_{(0,1)^{n-1}}
\widetilde{W}_{\Hom}(t, D_{\alpha}u+\overline F)\dx_{\alpha}:
\\ \nonumber&&\qquad
u_{|(0,1)^{n-1}}\in\wup_\loc((0,1)^{n-1};\R^m),\
u\ 1\hbox{-periodic\, in}\, x_{\alpha}\Bigr\}\Bigr)\dt\\
&\ge&
\int_0^1  \widetilde{W}_{\Hom}(t,\overline F)\dt
\end{eqnarray*}
by Jensen's inequality.
\qed

\begin{theorem}
For all $\overline F\in \M^{m\times (n-1)}$ we have 
$W_0(\overline F)=\overline W_{\Hom}(\overline F)$.
\end{theorem}

{\sc Proof.}
With fixed $\eta>0$ let $v$ be a test function for (\ref{a}) such that
$$
\int_{(0,1)^n}W_{\Hom}(x_n, Dv+\overline F)\dx\le
\overline W_{\Hom}(\overline F)+\eta\,.
$$
By Theorem \ref{oscihom} there
exists a sequence $v_j$ converging to $v$
such that $v_j=v$ on $\partial (0,1)^n$ (and, hence, in particular
$v_j$ is $1$-periodic in $x_\alpha$) and
\begin{eqnarray}\label{b}
\int_{(0,1)^n}W_{\Hom}(x_n, Dv+\overline F)\dx=\lim_{j\to +\infty}
\int_{\Omega_j\cap (0,1)^n}W(Dv_j+\overline F)\dx\,.
\end{eqnarray}
If we define
$u_j(x_{\alpha},x_n)={1\over j}v_j(jx_{\alpha},x_n)$ then $u_j\to 0$ in
$\lp((0,1)^n;\R^m)$ and
\begin{eqnarray}\label{d}
\nonumber
\int_{\Omega_j\cap (0,1)^n}W(Dv_j+\overline F)\dx&=&
{1\over j^{n-1}}\int_{\Omega_j\cap ((0,j)^{n-1}\times (0,1))}
W(Dv_j+\overline F)\dx\\
&=&
\nonumber
\int_{\Omega_{j^2}\cap (0,1)^n}W(Dv_j(jy_{\alpha},y_n)+\overline F)\dy\\
&=&
\nonumber
\int_{\Omega_{j^2}\cap (0,1)^n}W(D_{\alpha}u_j+\overline F|jD_nu_j)\dy\\
&=&
E_j(u_j+\overline F x_{\alpha},(0,1)^n);
\end{eqnarray}
hence, we can conclude that
\begin{eqnarray*}
W_0(\overline F)&\le& \liminf_{j\to +\infty}
E_j(u_j+\overline F x_{\alpha},(0,1)^n)\\
&=&
\liminf_{j\to +\infty}
\int_{\Omega_j\cap (0,1)^n}W(Dv_j+\overline F)\dx\\
&\le&
\overline W_{\Hom}(\overline F)+\eta
\end{eqnarray*}
by (\ref{inf}), (\ref{d}), (\ref{b}) and (\ref{a}).

\bigskip
Now we prove the converse inequality. Let $u_j\to 0$ be such that
$$
W_0(\overline F)=\lim_{j\to +\infty}E_j(u_j+\overline F x_{\alpha},(0,1)^n).
$$
By \cite{BFF} Lemma 2.6 we can choose $u_j$ $1$-periodic in $x_\alpha$;
let $v_j$ be defined by $v_j(x)=j u_j(x_\alpha/j,x_n)$.
With fixed $j,N\in \NN$,
$(0,1)^n=\bigcup_{m=1}^{N}(0,1)^{n-1}\times ((m-1)/N, m/N)$;
we can define a function $v_{j,m}$ by setting
$$
v_{j,m}(x_{\alpha},x_n)=\cases{
v_j(x_{\alpha},x_n+{2k\over N})
& if ${m-1\over N}- {2k\over N}< x_n<{m\over N}-{2k\over N} $
\cr\cr
v_j(x_{\alpha},{2m\over N}-x_n-{2k+2\over N}) & if
${m-1\over N}-{2k+1\over N} < x_n< {m\over N}-{2k+1\over N}$
}
$$
for $k\in \ZZ$, which is $1$-periodic in $x_{\alpha}$ and
$2/N$-periodic in $x_n$.
Hence, we can construct
$$
{w_{j,k}}_{|(0,1)^{n-1}\times ((m-1)/N, m/N)}=v_{j,m,k}(x)
$$
where $v_{j,m,k}(x)= {j\over k} v_{j,m}({k\over j}x)$,
such that $w_{j,k}$ is ${j\over k}$-periodic in $x_{\alpha}$ and
$$ {w_{j,k}}_{|(0,1)^{n-1}\times ((m-1)/N, m/N)} \rightarrow
\Bigl(0, \Bigl(\int_{(0,1)^n}D_n v_{j,m}dx\Bigr)x_n\Bigr)=w^m$$
as $k\to +\infty$, in $\lp ((0,1)^n;\R^m)$.
In this case the functions $w_{j,k}$ defined as above belong to
$\wup(\Omega_k\cap (0,1)^n;\R^m)$.

Finally, we define $w$ such that
$$
w_{|(0,1)^{n-1}\times ((m-1)/N, m/N)}=w^m$$
which is $1$-periodic in $x_{\alpha}$.
Let
$$A_j^{m/N}=\Omega_j\cap \{x_n=m/N\} $$
and
$$A_k^{m/N}=\Omega_k\cap \{x_n=m/N\},$$
we define
$$
E_j^N=\bigcup_{m=1}^N A_j^{m/N}\times ((m-1)/N, m/N)
$$
and
$$
E_k^N=\bigcup_{m=1}^N A_k^{m/N}\times ((m-1)/N, m/N).
$$
We restrict our analysis to the case where ${k/ j}$ odd,
the other case being dealt with by introducing a small error term.
Hence, if we use the notation
$$
I_l(u,(0,1)^n)=
\int_{E_l^N\cap (0,1)^n} W(Du )\dx
$$ 
($l=j$ or $k$) we have that
\begin{equation}\label{f2}
I_j(v_j+\overline Fx_{\alpha},(0,1)^n)=
I_k(w_{j,k}+\overline F x_{\alpha},(0,1)^n).
\end{equation}
Reasoning as in Theorems \ref{8a} and \ref{9} we get that
\begin{eqnarray*}
&&
I_{\Hom}(w+\overline F x_{\alpha},(0,1)^n)=
\Gamma\hbox{-}\lim_{k\to +\infty}
I_k(w+\overline Fx_{\alpha},(0,1)^n)\\
&=&
\sum_{m=1}^N \int_{(0,1)^{n-1}\times ((m-1)/N, m/N)}
W_{\Hom}(m/N, Dw+\overline F)\dx\\
&=&
\sum_{m=1}^N \int_{(0,1)^{n-1}\times ((m-1)/N, m/N)}
W_{\Hom}\Bigl({[x_n N]+1\over N}, Dw+\overline F\Bigr)\dx\\
&=&
\int_{(0,1)^n}W_{\Hom}\Bigl({[x_n N]+1\over N}, Dw+\overline F\Bigr)\dx\\
&\ge&
\int_0^1 \widetilde{W}_{\Hom}\Bigl({[x_n N]+1\over N},\overline F\Bigr)\dx_n
\end{eqnarray*}
by (\ref{c}).
Taking  the limit as $N\to +\infty$, we obtain
\begin{equation}\label{e}
I_{\Hom}(w+\overline F x_{\alpha},(0,1)^n)\ge \overline W_{\Hom}(\overline F)
\end{equation}
by Proposition \ref{gg}.
Hence,
\begin{eqnarray*}
E_j(u_j+\overline F x_{\alpha},(0,1)^n)
&=& \int_{\Omega_j\cap (0,1)^n}W(Dv_j+\overline F)\dx\\
&\ge&
\liminf_{k\to +\infty}I_{k}(w_{j,k}+\overline F x_{\alpha},(0,1)^n)\\
&\ge&
\overline W_{\Hom}(\overline F)
\end{eqnarray*}
by (\ref{d})-(\ref{e}).
By the choice of $(u_j)$ we get the desired inequality.
\qed

The proof of Theorem \ref{tfft} will be complete once we observe that
in the convex case formula (\ref{true}) simplifies that in Theorem \ref{9}
(see {\it e.g.} \cite{BD} Section 14.3).

\subsection{Convergence of minimum problems}
As an application of the $\Gamma$-convergence result of the previous section, we
describe the asymptotic behaviour of problems of the form
\begin{eqnarray}\label{med}\nonumber
m_{\e,\delta}&=&\min\Bigl\{\int_{\Omega(\e,\delta)} W(Du)\dx:
 u\in\lp(\omega\times(-\e,\e);\rr^m),\
\\
&&
u_{|\Omega(\e,\delta)}\in\wup(\Omega(\e,\delta);\rr^m),\
u=\phi\hbox{ on }(\partial\omega)\times(-\e,\e)
\Bigr\},
\end{eqnarray}
where $\phi=\phi(x_\alpha)\in\wup(\omega;\rr^m)$,
$\Omega(\e,\delta)$ is given by  (\ref{1.02}) and $f$ and $W$ satisfy the
hypotheses
of Theorem \ref{tfft}. By using Poincar\'e's inequality it can immediately
be checked that problem (\ref{med}) admits at least one solution for each
choice of $\e,\delta>0$.
The asymptotic behaviour of these solutions when $\e\to 0$ and $\delta<<\e$
is given
by the following result.

\begin{proposition}\label{mincon}
Let $\e$ and $\delta=\delta(\e)$ satisfy the hypotheses of Theorem
{\rm\ref{tfft}},
and for each $\e$ let $u_\e$ be a solution of {\rm(\ref{med})}. Then, upon
extracting
a subsequence, there exist a sequence $(v_\e)$ in
$\lp(\omega\times(-1,1);\rr^m)$
and a function $w\in\wup(\omega;\rr^m)$
such that

{\rm(i)} $v_\e= u_\e$ on $\Omega(\e,\delta(\e))$,

{\rm(ii)} if $w_\e(x_\alpha,x_n)=v_\e(x_\alpha,\e x_n)$, then $w_\e$ converges
(with the identification $w(x)=w(x_\alpha)$) to
$w$ in $\lp((\omega\times(-1,1);\rr^m)$,

{\rm(iii)} $w$ is a solution of the minimum problem
\begin{equation}\label{medo}
\widetilde m_{0}=\min\Bigl\{\int_{\omega} 2\overline W_\Hom(D_\alpha
u)\dx_\alpha:
u\in\lp(\omega;\rr^m),\
u=\phi\hbox{ on } \partial\omega
\Bigr\},
\end{equation}
where $\overline W_\Hom$ is defined by {\rm (\ref{treat}) and (\ref{true})},

{\rm(iv)}  $m_{\e,\delta(\e)}/\e$ converges to $\widetilde m_0$.
\end{proposition}

{\sc Proof.}
Note that, in the notation of Theorem \ref{tfft},
$\widetilde u_\e$ defined by  $\widetilde
u_\e(x_\alpha,x_n)=u_\e(x_\alpha,\e x_n)$ is a solution
of
\begin{eqnarray}\label{medr}\nonumber
\widetilde m_{\e}&=&{1\over\e}m_{\e,\delta(\e)}
=\min\Bigl\{\int_{\Omega_\e} W\Bigl(D_\alpha u,{1\over\e} D_n u\Bigr)\dx:
u\in\lp(\omega\times(-1,1);\rr^m),  \\
&&\qquad
u_{|\Omega_\e}\in\wup(\Omega_\e;\rr^m),\
u=\phi\hbox{ on }(\partial\omega)\times(-1,1)
\Bigr\}.
\end{eqnarray}
By \cite{BFF} Remark 2.3, upon extracting a subsequence, there exist
$w_\e\in\lp((\omega\times(-1,1);\rr^m)$ converging to some $w$ in
$\lp((\omega\times(-1,1);\rr^m)$, $D_nw=0$ and
$w_\e=\widetilde u_\e$ on $\Omega_\e$. By the well-known property of the
convergence of
minima and minimizers of $\Gamma$-converging functionals (see {\it e.g.}
\cite{BD} Theorem 7.2),
(iii) and (iv) follow from Theorem \ref{tfft}, since the $\Gamma$-limit is
not influenced by the boundary value $\phi$ (see \cite {BFF} Lemma 2.6).
\qed

\bigskip\noindent{\bf Acknowledgements}\ \
We gratefully acknowledge stimulating discussions with I.~Fon\-se\-ca,
and a very careful reading of the manuscript by the anonymous referee.
The research of AB was partially supported by Marie-Curie
fellowship ERBFMBICT972023 of the European Union program ``Training and
Mobility of Research\-ers'',
and benefitted from the hospitality of the
Max-Planck Institute for Mathematics in the Sciences, Leipzig (Germany).


\begin{thebibliography}{99}


\bibitem{ABP} {\rm G.  Anzellotti, S.  Baldo and D.  Percivale,}
{Dimension-reduction in variational problems, asymptotic
development in $\Gamma$-convergence and thin structures in elasticity,}
{\it Asymptotic Anal.} {\bf 9} (1994), 61--100.

\bibitem{BB} {\rm K. Bhattacharya and A. Braides},
Thin films with many small cracks.
{Preprint SISSA 1999}.

\bibitem{BJ}{\rm K. Bhattacharya and R.D. James}, 
A theory of thin films of martensitic materials with applications to
microactuators, 
{\it J. Mech. Phys. Solids} {\bf 47} (1999), 531-576.
 
\bibitem{Br} {\rm A. Braides,} {Homogenization of some almost
periodic functional,} {\it Rend. Accad. Naz. Sci. XL} {\bf 103}
(1985), 313--322.

\bibitem{B} {\rm A. Braides}, {\it $\Gamma$-convergence for Beginners}, 
Oxford University Press, Oxford, to appear.

\bibitem{BD} {\rm A. Braides and A. Defranceschi}, {\it
Homogenization of Multiple Integrals}, Oxford University Press,
Oxford, 1998.

\bibitem{BF} {\rm A. Braides and I. Fonseca},
Brittle thin films,
{\it Appl. Math. Optim.}, to appear.

\bibitem{BFF} {\rm A. Braides, I. Fonseca and G. Francfort,}
3D-2D asymptotic analysis for inhomogeneous thin films, to appear.

\bibitem{BC} {\rm R. Brizzi and J.P. Chalot,}
{Boundary homogenization and Neumann boundary value problems},
{\it Ric. Mat.} {\bf 46} (1997), 341-387.

\bibitem{BU} {\rm G. Buttazzo,}
{\it Semicontinuity, Relaxation and Integral Representation in the
Calculus of Variations}, Longman, Harlow, 1989.

\bibitem{Ca} {\rm D. Caillerie,} {Thin elastic and periodic plates},
{\it Math. Meth Appl. Sci.} {\bf 6} (1984), 159--191.

\bibitem{CV} C. Castaign and M. Valadier, {\it Convex Analysis and
Measurable Multifunctions}, Springer Verlag, Berlin, 1977.

\bibitem{DM} {\rm G. Dal Maso,} {\it An Introduction to
$\Gamma$-convergence},  Birkh\"auser, Boston, 1993.

\bibitem{DF} {\rm E. De Giorgi and T. Franzoni,} {Su un tipo di
convergenza  variazionale}, {\it Atti Accad. Naz. Lincei Rend. Cl.
Sci. Mat. Fis. Natur.} {\bf 58} (1975), 842--850.

\bibitem{EG} {\rm L.C. Evans and R.F. Gariepy,} {\it Measure
Theory and Fine Properties of Functions}, CRC Press, Ann Harbor, 1992.

\bibitem{FF} {\rm I. Fonseca and G. Francfort,}
{ 3D-2D asymptotic analysis of an optimal design problem for thin films},
{\it J. reine angew. Math.} {\bf 505} (1998), 173-202.

\bibitem{FMP} I. Fonseca, S. M\"uller  and P. Pedregal,
Analysis of concentration and oscillation effects generated by gradients.
{\it SIAM J. Math. Anal.} {\bf 29} (1998), 736--756.

\bibitem{FrM} G.A. Francfort and S. M\"uller, 
Combined effects of homogenization and singular perturbations in elasticity. 
{\it J. reine angew. Math.} {\bf 454} (1994), 1--35. 

\bibitem{LDR} {\rm H. Le Dret and A. Raoult,} {The nonlinear
membrane model as variational limit of nonlinear three-dimensional
elasticity}, {\it J. Math. Pures Appl.} {\bf 74} (1995), 549--578.

\bibitem{KV} {\rm R.V. Kohn and M. Vogelius,} {A new model for thin
plates with rapidly varying thickness. II: a convergence proof}, {\it
Quarterly Appl. Math.} {\bf 43} (1985), 1--22.

\bibitem{Ma} {\rm P. Marcellini}, {Periodic solutions and
homogenization of nonlinear  variational problems}, {\it Ann. Mat.
Pura Appl.} {\bf 117} (1978), 481--498.

\bibitem{Mu} {\rm S. M\"uller,} {Homogenization of nonconvex
integral functionals and cellular elastic materials}, {\it Arch.
Rational Mech. Anal.} {\bf 99} (1987), 189--212.

\bibitem{Shu} {\rm Y.C. Shu,} {Heterogeneous thin films of martensitic
materials}, Preprint, 1998.


\end{thebibliography}
\end{document}